\documentclass[smallcondensed]{svjour3}
\smartqed
\usepackage{graphicx}
\usepackage{epsfig,multirow}
\usepackage{amsmath,amssymb,bm}
\usepackage{amsfonts,dsfont,color,bbm}
\usepackage{algorithm,algorithmic,multicol}
\usepackage{epsf,epstopdf}

\usepackage{natbib}
\usepackage{subcaption}
\usepackage{latexsym,amssymb,amsmath,graphicx,epsf,bbm,bm}
\usepackage{fixltx2e} % for figures to appear in  correct order
\usepackage{relsize}
\usepackage[font=footnotesize,labelfont=bf]{caption}
\usepackage{slashbox}
\usepackage{epsfig}
\usepackage{graphicx}
\usepackage{algorithm,algorithmic}
\usepackage{epstopdf,epsf,epsfig}

\usepackage{mathtools}

\bibpunct{(}{)}{;}{a}{}{,}

\renewcommand{\vec}[1]{\mbox{\boldmath${#1}$}}

\newcommand{\ei}{\end{itemize}}
\newcommand{\bi}{\begin{itemize}}

	\newcommand{\vI}{\vec{I}}
	
	\newcommand{\vq}{\vec{q}}

	\newcommand{\vw}{\vec{w}}

	\newcommand{\vz}{\vec{z}}

	\newcommand{\vy}{\vec{y}}

	\newcommand{\vx}{\mbox{$\vec{x}$}}

	\newcommand{\bR}{\mathbf{R}}

	\newcommand{\MB}{\left[\begin{array}}
		\newcommand{\ME}{\end{array}\right]}

	\newcommand{\nunder}[2][5]{\mathrlap{\mkern\the\numexpr#1/2mu\relax\underline{\phantom{\mathrm{#2}\mkern-#1mu}}}#2}
	
	\newcommand{\Real}{\mathbbm{R}}
	\newcommand{\lmtk}{\lambda_t^{(k)}}
	\newcommand{\ltk}{\l_t^{(k)}}
	\newcommand{\dltk}{\delta_{t-1}^{(k)}}

\begin{document}

\title{A Novel Family of Boosted Online Regression Algorithms with Strong Theoretical Bounds}

\titlerunning{Boosted Online Regression Algorithms with Strong Theoretical Bounds}        % if too long for running head

\author{Dariush Kari \and Farhan Khan \and Selami Ciftci \and Suleyman S. Kozat}

\institute{Dariush Kari \and Suleyman S. Kozat \at
              Department of Electrical and Electronics Engineering, Bilkent University \\
              Ankara 06800, Turkey \\
              \email{kari@ee.bilkent.edu.tr, kozat@ee.bilkent.edu.tr}\\
              \\
			Selami Ciftci\at
			Turk Telekom Communications Services Inc., Istanbul, Turkey\\
			\email{selami.ciftci1@turktelekom.com.tr}\\
			\\
			Farhan Khan \at
			Department of Electrical and Electronics Engineering, Bilkent University \\
			Ankara 06800, Turkey \\
			\email{khan@ee.bilkent.edu.tr}\\
			 and also Electrical Engineering Department, COMSATS Institute of Information Technology, Pakistan \\
	\email{engrfarhan@ciit.net.pk}
}

% \date{Received: date / Accepted: date}

\maketitle

\begin{abstract}
	We investigate boosted online regression and propose a novel family of regression algorithms with strong theoretical bounds. In addition, we implement several variants of the proposed generic algorithm. We specifically provide theoretical bounds for the performance of our proposed algorithms that hold in a strong mathematical sense. We achieve guaranteed performance improvement over the conventional online regression methods without any statistical assumptions on the desired data or feature vectors. We demonstrate an intrinsic relationship, in terms of boosting, between the adaptive mixture-of-experts and data reuse algorithms. Furthermore, we introduce a boosting algorithm based on random updates that is significantly faster than the conventional boosting methods and other variants of our proposed algorithms while achieving an enhanced performance gain. Hence, the random updates method is specifically applicable to the fast and high dimensional streaming data. Specifically, we investigate Newton Method-based and Stochastic Gradient Descent-based linear regression algorithms in a mixture-of-experts setting, and provide several variants of these well known adaptation methods. However, the proposed algorithms can be extended to other base learners, e.g., nonlinear, tree-based piecewise linear. Furthermore, we provide theoretical bounds for the computational complexity of our proposed algorithms. We demonstrate substantial performance gains in terms of mean square error over the base learners through an extensive set of benchmark real data sets and simulated examples.

\end{abstract}

\section{Introduction}\label{sec:introduction}
Boosting is considered as one of the most important ensemble learning methods in the machine learning literature and it is extensively used in several different real life applications from classification to regression (\cite{Bauer1999, Dietterich2000, Schapire1999, freund_book, adaboost, regression_boost, Shalev-Shwartz2010, Saigo2009, Demiriz2002}). As an ensemble learning method (\cite{Fern2003, Soltanmohammadi2016, duda_book}), boosting combines several parallel running ``weakly'' performing algorithms to build a final ``strongly'' performing algorithm (\cite{Soltanmohammadi2016, Freund2001, freund_book, Mannor2002}). This is accomplished by finding a linear combination of weak learning algorithms in order to minimize the total loss over a set of training data commonly using a functional gradient descent (\cite{Duffy2002, adaboost}). Boosting is successfully applied to several different problems in the machine learning literature including classification (\cite {Jin2007, Chapelle2011, adaboost}), regression (\cite{Duffy2002, regression_boost}), and prediction (\cite{pred_boost, pred_boost2}). However, significantly less attention is given to the idea of boosting in online regression framework. To this end, our goal is (a) to introduce a new boosting approach for online regression, (b) derive several different online regression algorithms based on the boosting approach, (c) provide mathematical guarantees for the performance improvements of our algorithms, and (d) demonstrate the intrinsic connections of boosting with the adaptive mixture-of-experts algorithms (\cite{adamix, kozat}) and data reuse algorithms (\cite{datareuse1}).\par
Although boosting is initially introduced in the batch setting (\cite{adaboost}), where algorithms boost themselves over a fixed set of training data, it is later extended to the online setting (\cite{ozaboost}). In the online setting, however, we neither need nor have access to a fixed set of training data, since the data samples arrive one by one as a stream (\cite{Ben-David1997,Fern2003,Lu2016}). Each newly arriving data sample is processed and then discarded without any storing. The online setting is naturally motivated by many real life applications especially for the ones involving big data, where there may not be enough storage space available or the constraints of the problem require instant processing (\cite{big_data}). Therefore, we concentrate on the online boosting framework and propose several algorithms for online regression tasks. In addition, since our algorithms are online, they can be directly used in adaptive filtering applications to improve the performance of conventional mixture-of-experts methods (\cite{adamix}). For adaptive filtering purposes, the online setting is especially important, where the sequentially arriving data is used to adjust the internal parameters of the filter, either to dynamically learn the underlying model or to track the nonstationary data statistics (\cite{adamix,sayed_book}).\par
Specifically, we have $m$ parallel running weak learners (WL) (\cite{freund_book}) that receive the input vectors sequentially. Each WL uses an update method, such as the second order Newton's Method (NM) or Stochastic Gradient Descent (SGD), depending on the target of the applications or problem constraints (\cite{sayed_book}). After receiving the input vector, each algorithm produces its output and then calculates its instantaneous error after the observation is revealed. In the most generic setting, this estimation/prediction error and the corresponding input vector are then used to update the internal parameters of the algorithm to minimize a priori defined loss function, e.g., instantaneous error for the SGD algorithm. These updates are performed for all of the $m$ WLs in the mixture. However, in the online boosting approaches, these adaptations at each time proceed in rounds from top to bottom, starting  from the first WL to the last one to achieve the ``boosting'' effect (\cite{onlineboost}). Furthermore, unlike the usual mixture approaches (\cite{adamix,kozat}), the update of each WL depends on the previous WLs in the mixture. In particular, at each time $t$, after the $k^{th}$ WL calculates its error over $(\vec{x}_t,d_t)$ pair, it passes a certain weight to the next WL, the $(k+1)^{th}$ WL, quantifying how much error the constituent WLs from $1^{st}$ to $k^{th}$ made on the current $(\vec{x}_t,d_t)$ pair. Based on the performance of the WLs from $1$ to $k$ on the current $(\vec{x}_t,d_t)$ pair,  the $(k+1)^{th}$ WL may give a different emphasis (importance weight) to $(\vec{x}_t,d_t)$ pair in its adaptation in order to rectify the mistake of the previous WLs.\par
The proposed idea for online boosting is clearly related to the adaptive mixture-of-experts algorithms widely used in the machine learning literature, where several parallel running adaptive algorithms are combined to improve the performance. In the mixture methods, the performance improvement is achieved due to the diversity provided by using several different adaptive algorithms each having a different view or advantage (\cite{kozat}). This diversity is exploited to yield a final combined algorithm, which achieves a performance better than any of the algorithms in the mixture. Although the online boosting approach is similar to mixture approaches (\cite{kozat}), there are significant differences. In the online boosting notion, the parallel running algorithms are not independent, i.e., one deliberately introduces the diversity by updating the WLs one by one from the first WL to the $m^{th}$ WL for each new sample based on the performance of all the previous WLs on this sample. In this sense, each adaptive algorithm, say the $(k+1)^{th}$ WL, receives feedback from the previous WLs, i.e., $1^{st}$ to $k^{th}$, and updates its inner parameters accordingly. As an example, if the current $(\vec{x}_t,d_t)$ is well modeled by the previous WLs, then the $(k+1)^{th}$ WL performs minor update using $(\vec{x}_t,d_t)$ and may give more emphasis (importance weight) to the later arriving  samples that may be worse modeled by the previous WLs. Thus, by boosting, each adaptive algorithm in the mixture can concentrate on different parts of the input and output pairs achieving diversity and significantly improving the gain.\par 
%Moreover, we introduce the ``random updates'' method for boosting, which significantly reduces the computational complexity, while achieving the performance of other variants of our boosting method. This is because in this scenario, the $k^{th}$ WL is updated with probability $\lmtk$, which depends on the performance of other filters.\par
The linear online learning algorithms, such as SGD or NM, are among the simplest as well as the most widely used regression algorithms in the real-life applications (\cite{sayed_book}). Therefore, we use such algorithms as base WLs in our boosting algorithms. To this end, we first apply the boosting notion to several parallel running linear NM-based WLs and introduce three different approaches to use the importance weights (\cite{onlineboost}), namely ``weighted updates'',``data reuse'', and ``random updates''. In the first approach, we use the importance weights directly to produce certain weighted NM algorithms. In the second approach, we use the importance weights to construct data reuse adaptive algorithms (\cite{ozaboost}). However, data reuse in boosting, such as (\cite{ozaboost}), is significantly different from the usual data reusing approaches in adaptive filtering (\cite{datareuse1}). As an example, in boosting, the importance weight coming from the $k^{th}$ WL determines the data reuse amount in the $(k+1)^{th}$ WL, i.e., it is not used for the $k^{th}$ filter, hence, achieving the diversity. The third approach uses the importance weights to decide whether to update the constituent WLs or not, based on a random number generated from a Bernoulli distribution with the parameter equal to the weight. The latter method can be effectively used for big data processing (\cite{bigdata}) due to the reduced complexity. The output of the constituent WLs is also combined using a linear mixture algorithm to construct the final output. We then update the final combination algorithm using the SGD algorithm (\cite{kozat}). Furthermore, we extend the boosting idea to parallel running linear SGD-based algorithm similar to the NM case.\par
We start our discussions by investigating the related works in Section \ref{sec:rel_work}. We then introduce the problem setup and background in Section \ref{sec:problem}, where we provide individual sequence as well as MSE convergence results for the NM and SGD algorithms. We introduce our generic boosted online regression algorithm in Section \ref{sec:BORA} and provide the mathematical justifications for its performance. Then, in Sections \ref{sec:BRLS} and \ref{sec:BLMS}, three different variants of the proposed boosting algorithm are derived, using the NM and SGD, respectively. Then, in Section \ref{sec:analysis} we provide the mathematical analysis for the computational complexity of the proposed algorithms. The paper concludes with extensive sets of experiments over the well known benchmark data sets and simulation models widely used in the machine learning literature to demonstrate the significant gains achieved by the boosting notion.

\section{Related Works}\label{sec:rel_work}
AdaBoost is one of the earliest and most popular boosting methods, which has been used for binary and multiclass classifications as well as regression (\cite{adaboost}). This algorithm has been well studied and has clear theoretical guarantees, and its excellent performance is explained rigorously (\cite{breiman97}). However, AdaBoost cannot perform well on the noisy data sets (\cite{smooth_boost}), therefore, other boosting methods have been suggested that are more robust against noise.\par
In order to reduce the effect of noise, SmoothBoost was introduced in (\cite{smooth_boost}) in a batch setting. Moreover, in (\cite{smooth_boost}) the author proves the termination time of the SmoothBoost algorithm by simultaneously obtaining upper and lower bounds on the weighted advantage of all samples over all of the weak learners. We note that the SmoothBoost algorithm avoids overemphasizing the noisy samples, hence, provides robustness against noise. In (\cite{ozaboost}), the authors extend bagging and boosting methods to an online setting, where they use a Poisson sampling process to approximate the reweighting algorithm. However, the online boosting method in (\cite{ozaboost}) corresponds to AdaBoost, which is susceptible to noise. In (\cite{babenko}), the authors use a greedy optimization approach to develop the boosting notion to the online setting and introduce stochastic boosting. Nevertheless, while most of the online boosting algorithms in the literature seek to approximate AdaBoost, (\cite{onlineboost}) investigates the inherent difference between batch and online learning, extend the SmoothBoost algorithm to an online setting, and provide the mathematical guarantees for their algorithm. (\cite{onlineboost}) points out that the online weak learners do not need to perform well on all possible distributions of data, instead, they have to perform well only with respect to smoother distributions. Recently, in (\cite{beygelzimer15classification}) the authors have developed two online boosting algorithms for classification, an optimal algorithm in terms of the number of weak learners, and also an adaptive algorithm using the potential functions and boost-by-majority (\cite{boost_majority}).\par
In addition to the classification task, the boosting approach has also been developed for the regression (\cite{Duffy2002}). In (\cite{bertoni97}), a boosting algorithm for regression is proposed, which is an extension of Adaboost.R (\cite{bertoni97}). Moreover, in (\cite{Duffy2002}), several gradient descent algorithms are presented, and some bounds on their performances are provided. In (\cite{babenko}) the authors present a family of boosting algorithms for online regression through greedy minimization of a loss function. Also, in (\cite{beygelzimer15regression}) the authors propose an online gradient boosting algorithm for regression.\par
In this paper we propose a novel family of boosted online algorithms for the regression task using the ``online boosting'' notion introduced in (\cite{onlineboost}), and investigate three different variants of the introduced algorithm. Furthermore, we show that our algorithm can achieve a desired mean squared error (MSE), given a sufficient amount of data and a sufficient number of weak learners. In addition, we use similar techniques to (\cite{smooth_boost}) to prove the correctness of our algorithm. We emphasize that our algorithm has a guaranteed performance in an individual sequence manner, i.e., without any statistical assumptions on the data. In establishing our algorithm and its justifications, we refrain from changing the regression problem to the classification problem, unlike the AdaBoost.R (\cite{adaboost}). Furthermore, unlike the online SmoothBoost (\cite{onlineboost}), our algorithm can learn the guaranteed MSE of the weak learners, which in turn improves its adaptivity.\par

%For this case, we combine the outputs of the constituent filters using a linear filter, which is trained using the LMS algorithm. 

%However, note that although linear filters have a low complexity, piecewise linear filters deliver a significantly superior performance in real life applications \cite{denizcan}, with a comparable complexity. These filters mitigate the overfitting, stability and convergence issues tied to nonlinear models \cite{sayed_book,Singer3, Nascimento1, sayed2}, while effectively improving the modeling power relative to linear filters \cite{denizcan}.

%For all these different cases, we derive the corresponding mean squared error (MSE) results and provide performance bounds in an individual sequence manner \cite{sinkoz, cesa_book}.\par

\section{Problem Description and Background} \label{sec:problem}

All vectors are column vectors and represented by bold lower case letters. Matrices are represented by bold upper case letters. For a vector $\vec{a}$ (or a matrix $\vec{A}$),  $\vec{a}^T$ (or $\vec{A}^T$) is the transpose and Tr($\vec{A}$) is the trace of the matrix $\vec{A}$. Here, $\vec{I}_m$ and $\vec{0}_m$ represent the identity matrix of dimension $m\times m$ and the all zeros vector of length $m$, respectively. Except $\vec{I}_m$ and $\vec{0}_m$, the time index is given in the subscript, i.e., $x_t$ is the sample at time $t$. We work with real data for notational simplicity. We denote the mean of a random variable $x$ as $E[x]$. Also, we show the cardinality of a set $S$ by $|S|$.

We sequentially receive $r$-dimensional input (regressor) vectors $\{\vec{x}_t\}_{t \geq 1}$, $\vec{x}_t \in {\Real}^r$, and desired data $\{d_t\}_{t \geq 1}$, and estimate $d_t$ by $\hat{d}_t = f_t(\vx_t)$, where $f_t(.)$ is an online regression algorithm. At each time $t$ the estimation error is given by $e_t=d_t-\hat{d}_t$ and is used to update the parameters of the WL. For presentation purposes, we assume that $d_t \in [-1,1]$, however, our derivations hold for any bounded but arbitrary desired data sequences. In our framework, we do not use any statistical assumptions on the input feature vectors or on the desired data such that our results are guaranteed to hold in an individual sequence manner (\cite{Ko08}).\par 
%However, we also provide steady-state, tracking and transient MSE analysis of the algorithms under widely used statistical models in the signal processing literature \cite{sayed_book} for completeness.\par
%Note that although nonlinear filters can outperform linear filters, they usually undergo overfitting, stability, and convergence issues \cite{denizcan,kozat2}. Furthermore, nonlinear filters generally have higher computational complexities, which limits their use in most of the real-life applications \cite{denizcan,kozat2}. To overcome these problems, piecewise linear filters are proposed, which mitigate the overfitting and stability issues, while offering a comparable modeling performance to the nonlinear filters \cite{denizcan,kozat2}. Therefore, in this paper, we are particularly interested in piecewise linear filters, which serve as an elegant alternative to linear filters. Nevertheless, for illustration, we first explain the basic principles of linear filters, and their extension to the piecewise linear case. Then, in Sections \ref{sec:BRLS} and \ref{sec:BLMS}, we introduce our algorithm in a piecewise linear model.\par
The linear methods are considered as the simplest online modeling or learning algorithms, which estimate the desired data $d_t$ by a linear model as $\hat{d}_t = \vw_t^{T} \vx_t$, where $\vw_t$ is the linear algorithm's coefficients at time $t$. Note that the previous expression also covers the affine model if one includes a constant term in $\vec{x}_t$, hence we use the purely linear form for notational simplicity. When the true $d_t$ is revealed, the algorithm updates its coefficients $\vec{w}_t$ based on the error $e_t$. As an example, in the basic implementation of the NM algorithm, the coefficients are selected to minimize the accumulated squared regression error up to time $t-1$ as
\begin{align}
\vec{w}_t & = \arg \min_{\vw} \sum_{l=1}^{t-1} (d_l - \vec{x}_l^T \vw)^2, \nonumber\\
& = \left( \sum_{l=1}^{t-1} \vec{x}_l \vec{x}_l^T \right)^{-1} \left( \sum_{l=1}^{t-1} \vec{x}_l d_l \right), \label{eq:rls}
\end{align}
where $\vw$ is a fixed vector of coefficients. The NM algorithm is shown to enjoy several optimality properties under different statistical settings (\cite{sayed_book}).  Apart from these results and more related to the framework of this paper, the NM algorithm is also shown to be rate optimal in an individual sequence manner (\cite{universal}). As shown in (\cite{universal}) (Section V), when applied to any sequence $\{\vec{x}_t\}_{t \geq 1}$ and $\{d_t\}_{t \geq 1}$, the accumulated squared error of the NM algorithm is as small as the accumulated squared error of the best batch least squares (LS) method that is directly optimized for these realizations of the sequences, i.e., for all $T$, $\{\vec{x}_t\}_{t \geq 1}$ and $\{d_t\}_{t \geq 1}$, the NM achieves
\begin{equation}
\sum_{l=1}^{T} (d_l - \vec{x}_l^T \vec{w}_l)^2 - \min_{\vw} \sum_{l=1}^{T} (d_l - \vec{x}_l^T \vw)^2 \leq O(\ln T). \label{eq:rls_bound}
\end{equation}
The NM algorithm is a member of the Follow-the-Leader type algorithms (\cite{cesa_book}) (Section 3), where one uses the best performing linear model up to time $t-1$ to predict $d_t$.  Hence, \eqref{eq:rls_bound} follows by direct application of the online convex optimization results (\cite{olsurvey}) after  regularization. The convergence rate (or the rate of the regret)  of the NM algorithm is also shown to be optimal so that $O(\ln T)$ in the upper bound cannot be improved (\cite{sinkoz}). It is also shown in (\cite{sinkoz}) that one can reach the optimal upper bound (with exact scaling terms) by using a slightly modified version of \eqref{eq:rls}
\begin{equation}
\vec{w}_t = \left( \sum_{l=1}^{t} \vec{x}_l \vec{x}_l^T \right)^{-1}\left( \sum_{l=1}^{t-1} \vec{x}_l d_l \right). \label{eq:mine}
\end{equation}
Note that the extension \eqref{eq:mine} of \eqref{eq:rls} is a forward algorithm (Section 5 of \cite{density}) and one can show that, in the scalar case, the predictions of \eqref{eq:mine} are always bounded (which is not the case for \eqref{eq:rls}) (\cite{sinkoz}).\par
We emphasize that in the basic application of the NM algorithm, all data pairs $(d_l, \vec{x}_l)$, $l=1,\ldots,t$, receive the same ``importance'' or weight in \eqref{eq:rls}. Although there exists exponentially weighted or windowed versions of the basic NM algorithm (\cite{sayed_book}), these methods weight (or concentrate on) the most recent samples for better modeling of the nonstationarity (\cite{sayed_book}). However, in the boosting framework (\cite{adaboost}), each sample pair receives a different weight based on not only those weighting schemes, but also the performance of the boosted algorithms on this pair. As an example, if a WL performs worse on a sample, the next WL concentrates more on this example to better rectify this mistake. In the following sections, we use this notion to derive different boosted online regression algorithms.\par
Although in this paper we use linear WLs for the sake of notational simplicity, one can readily extend our approach to nonlinear and piecewise linear regression methods. For example, one can use tree based online regression methods (\cite{farhan, Denizcan, kozat_ctw}) as the weak learners, and boost them with the proposed approach.
\vspace{-0.25cm}
\section{New Boosted Online Regression Algorithm}\label{sec:BORA}
In this section we present the generic form of our proposed algorithms and provide the guaranteed performance bounds for that. Regarding the notion of ``online boosting'' introduced in (\cite{onlineboost}), the online weak learners need to perform well only over smooth distributions of data points. We first present the generic algorithm in Algorithm (\ref{alg:BORA}) and provide its theoretical justifications, then discuss about its structure and the intuition behind it.
%Then, in Sections \ref{sec:BRLS} and \ref{sec:BLMS}, we respectively present the boosting algorithm for the SGD and the second order Newton methods with three variants for each.\par

\begin{algorithm}[htb]
	\caption{Boosted online regression algorithm}\label{alg:BORA}
	\begin{algorithmic}[1]
		\STATE Input: $(\vec{x}_t,d_t)$ (data stream), $m$ (number of weak learners running in parallel), $\sigma_m^2$ (the modified desired MSE), and $\sigma^2$ (the guaranteed achievable weighted MSE).
		\STATE Initialize the regression coefficients $\vec{w}_{1}^{(k)}$ for each WL; and the combination coefficients as $\vec{z}_1=\frac{1}{m}[1,1,\ldots,1]^T$;
		\FOR{$t=1$ \TO $T$}
		\STATE Receive the regressor data instance $\vec{x}_t$;
		\STATE Compute the WLs outputs $\hat{d}_t^{(k)}$;
		\STATE Produce the final estimate $\hat{d}_t = \vec{z}_t^T\vec{y}_{t} = \vec{z}_t^T[\hat{d}_t^{(1)}, \ldots, \hat{d}_t^{(m)}]^T$;
		\STATE Receive the true output $d_t$ (desired data);
		\STATE $\lambda_t^{(1)}=1$; $l_t^{(1)}=0$;
		\FOR{$k=1$ \TO $m$}
		\STATE $\lambda_t^{(k)} = \min \left\{ 1, \left(\sigma^2\right)^{\; l_t^{(k)}/2} \right\}$;
		\STATE Update the WL$^{(k)}$, such that it has a weighted MSE $\leq \sigma^2$;
		\STATE $e_t^{(k)}= d_t -\hat{d}_t^{(k)}$; 
		
		%\STATE $\delta_t^{(k)} = \frac{\Lambda_{t-1}^{(k)}\delta_{t-1}^{(k)}+\frac{\lambda_{t}^{(k)}}{4} \left(d_{t} -\left[f_{t}^{(k)}(\vx_{t})\right]^{+}\right)^2}{\Lambda_{t-1}^{(k)}+\lambda_{t}^{(k)}}$;
		%\STATE $\Lambda_{t}^{(k)} = \Lambda_{t-1}^{(k)}+\lambda_{t}^{(k)}$
		\STATE $l_t^{(k+1)} =  l_t^{(k)}+\left[ \sigma^2_m - \left(e_t^{(k)}\right)^2 \right]$;
		\ENDFOR
		\STATE Update $\vz_t$ based on $e_t=d_t-\vec{z}_t^T\vec{y}_{t}$;
		%\STATE $\vec{g}_t = \frac{\vec{P}_t \vec{y}_t}{\lambda+ \vec{y}_t^T \vec{P}_{t} \vec{y}_t }$;
		%\STATE $\vec{z}_{t+1} = \vec{z}_t + e_t \vec{g}_t$;
		%\STATE $\vec{P}_{t+1} = \lambda^{-1}\vec{P}_t-\lambda^{-1} \vec{g}_t \vec{y}_t^T \vec{P}_t$;
		\ENDFOR
	\end{algorithmic}
\end{algorithm}

In Algorithm \ref{alg:BORA}, we have $m$ copies of an online WL, each of which is guaranteed to have a weighted MSE of at most $\sigma^2$. We prove that the Algorithm 1 can reach a desired MSE, $\sigma_d^2$, through Lemma 1, Lemma 2, and Theorem 1. Note that since we assume $d_t \in [-1,1]$, the trivial solution $\hat{d}_t=0$ incurs an MSE of at most $1$. Therefore, we define a weak learner as an algorithm which has an MSE less than $1$.\par
\noindent
\textbf{Lemma 1. } {\em In Algorithm \ref{alg:BORA}, if there is an integer $M$ such that $\sum_{t=1}^{T} \lambda_t^{(k)} \geq \kappa T$ for every $k \leq M$, and also $\sum_{t=1}^{T} \lambda_t^{(M+1)} < \kappa T$, where $0<\kappa < \sigma_d^2$ is arbitrarily chosen, it can reach a desired MSE, $\sigma_d^2$.}\\
\noindent
\textbf{Proof. } The proof of Lemma 1 is given in Appendix \ref{app:lem1}.\par
\noindent
\textbf{Lemma 2. } {\em If the weak learners are guaranteed to have a weighted MSE less than $\sigma^2$, i.e.,
	\[
	\forall k: \qquad 
	\frac{\sum_{t=1}^T \lambda_t^{(k)}(e_t^{(k)})^2}{4\sum_{t=1}^T \lambda_t^{(k)}} \leq \sigma^2 \leq \frac{1}{4},
	\]
	there is an integer $M$ that satisfies the conditions in Lemma 1.}\par
\noindent
\textbf{Proof. } The proof of Lemma 2 is given in Appendix \ref{app:lem2}.\par
\noindent
\textbf{Theorem 1. }{\em If the weak learners in line 11 of Algorithm \ref{alg:BORA} achieve a weighted MSE of at most $\sigma^2 < \frac{1}{4}$ , there exists an upper bound for $m$ such that the algorithm reaches the desired MSE.}\\
\noindent
\textbf{Proof. } This theorem is a direct consequence of combining Lemma 1 and Lemma 2.\par
Note that although we are using copies of a base learner as the weak learners and seek to improve its performance, the constituent WLs can be different. However, by using the boosting approach, we can improve the MSE performance of the overall system as long as the WLs can provide a weighted MSE of at most $\sigma^2$. For example, we can improve the performance of mixture-of-experts algorithms (\cite{adamix}) by leveraging the boosting approach introduced in this paper.\par
\begin{figure}[!htb]
	\centering
	\includegraphics[width=0.8\textwidth]{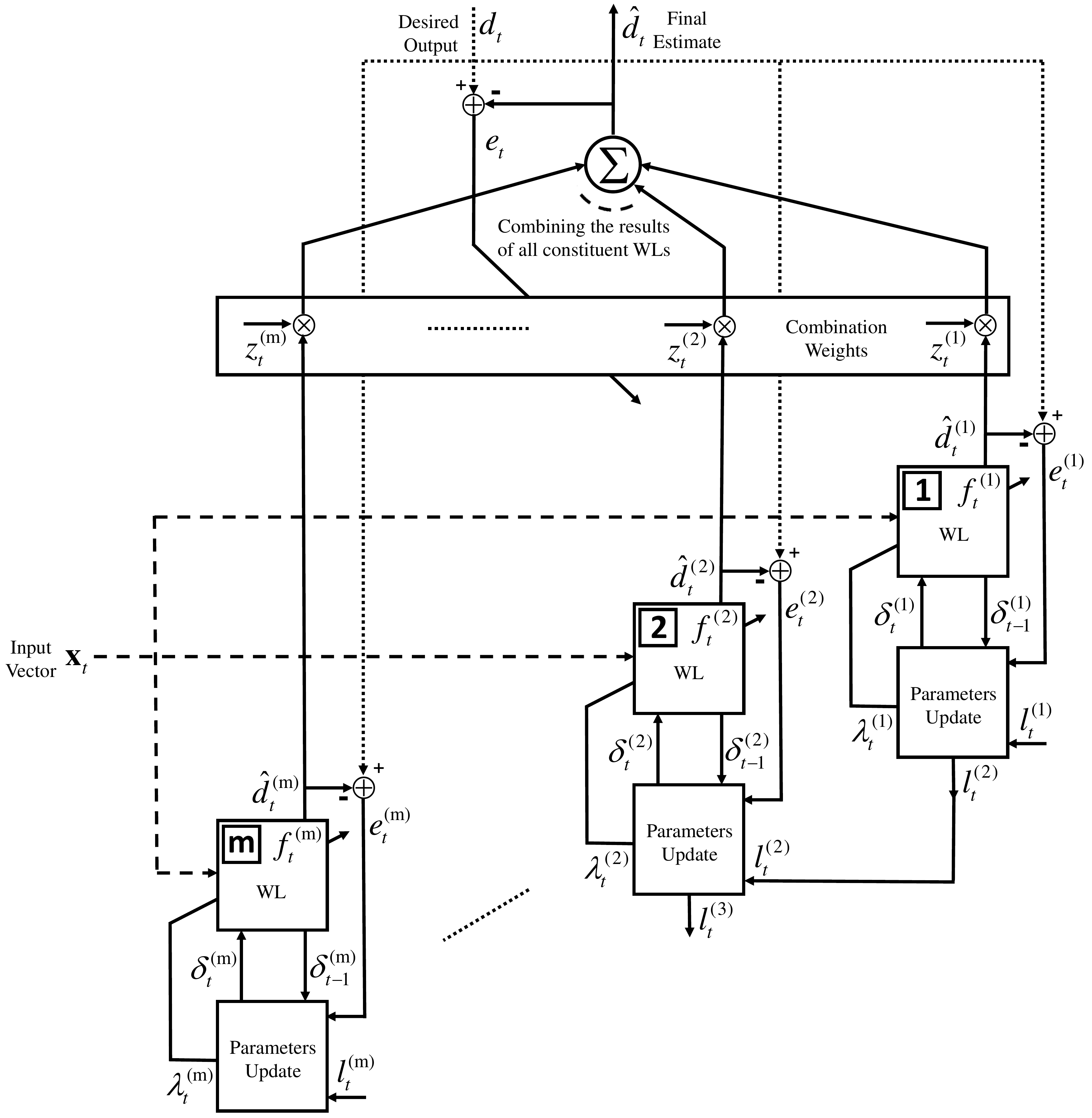}
	\caption{The block diagram of a boosted online regression system that uses the input vector $\vec{x}_t$ to produce the final estimate $\hat{d}_t$. There are $m$ constituent WLs $f_t^{(1)}, \ldots,f_t^{(m)}$, each of which is an online linear algorithm that generates its own estimate $\hat{d}_t^{(k)}$. The final estimate $\hat{d}_t$ is a linear combination of the estimates generated by all these constituent WLs, with the combination weights $z_t^{(k)}$'s corresponding to $\hat{d}_t^{(k)}$'s. The combination weights are stored in a vector which is updated after each iteration $t$. At time $t$ the $k^{th}$ WL is updated based on the values of $\lambda_t^{(k)}$ and $e_t^{(k)}$, and provides the $(k+1)^{th}$ filter with $l_t^{(k+1)}$ that is used to compute $\lambda_t^{(k+1)}$. The parameter $\delta_t^{(k)}$ indicates the weighted MSE of the $k^{th}$ WL over the first $t$ estimations, and is used in computing $\lambda_t^{(k)}$.} \label{fig:boost_block}
\end{figure}

As shown in Fig.~\ref{fig:boost_block}, at each iteration $t$, we have $m$ parallel running WLs with estimating functions $f_{t}^{(k)}$, producing estimates $\hat{d}_t^{(k)}= f_t^{(k)}(\vec{x}_t)$ of $d_t$, $k=1,\ldots,m$. As an example, if we use $m$ ``linear'' algorithms, $\hat{d}_t^{(k)} = \vec{x}_t^{T}\vw_t^{(k)}$ is the estimate generated by the $k^{th}$ WL. The outputs of these $m$ WLs are then combined using the linear weights $\vz_t$ to produce the final estimate as $\hat{d}_t = \vz_t^T\vec{y}_t$ (\cite{kozat}), where $\vec{y}_t \triangleq [\hat{d}_t^{(1)}, \ldots, \hat{d}_t^{(m)}]^T$ is the vector of outputs. After the desired output $d_t$ is revealed, the $m$ parallel running WLs will be updated for the next iteration. Moreover, the linear combination coefficients $\vec{z}_t$ are also updated using the normalized SGD (\cite{sayed_book}), as detailed later in Section~\ref{sec:final_rls}.\par

%\begin{figure}[!htb]
%	\centering
%	\includegraphics[width=0.3\textwidth]{figures/boost_detail.pdf}\\
%	\caption{Parameters update block of the $k$th constituent filter, which is embedded in the $k$th filter block as depicted in Fig. \ref{fig:boost_block}. This block receives the parameter $l_t^{(k)}$ provided by the $(k-1)$th filter, and uses that in computing $\lambda_t^{(k)}$. It also computes $l_t^{(k+1)}$ and passes it to the $(k+1)$th filter. The parameter $[e_t^{(k)}]^{+}$ represents the error of the thresholded estimate as explained in \eqref{eq:delta}, and $\Lambda_t^{(k)}$ shows the sum of the weights $\lambda_1^{(k)}, \ldots,\lambda_t^{(k)}$. The WMSE parameter $\delta_{t-1}^{(k)}$ represents the time averaged weighted square error made by the $k$th filter up to time $t-1$.}\label{fig:boost_detail}
%\end{figure}
%
%\begin{figure}[!htb]
%	\centering
%	\includegraphics[width=0.3\textwidth]{figures/piecewise.pdf}\\
%	\caption{A sample piecewise linear adaptive filter, used as the $k$th constituent filter in the system depicted in Fig. \ref{fig:boost_block}. This fliter consists of $N$ linear filters, one of which produces the estimate at each iteration $t$. Based on where the input vector at time $t$, $\vec{x}_t$, lies in the input vector space, one of the $s_{i,t}^{(k)}$'s is 1 and all others are 0. Hence, at each iteration only one of the linear filters is used for estimation and upadated correspondingly.}\label{fig:piece}
%\end{figure}

After $d_t$ is revealed, the constituent WLs, $f_t^{(k)}$, $k=1,\ldots,m$, are consecutively updated, as shown in Fig.~\ref{fig:boost_block}, from top to bottom, i.e., first $k=1$ is updated, then, $k=2$ and finally $k=m$ is updated. However, to enhance the performance, we use a boosted updating approach (\cite{adaboost}), such that the $(k+1)^{th}$ WL receives a ``total loss'' parameter, $l_t^{(k+1)}$, from the $k^{th}$ WL, as
\begin{equation}
l_t^{(k+1)} =  l_t^{(k)}+\left[ \sigma_m^2 - \left(d_t -f_t^{(k)}(\vx_t)\right)^2 \right], \label{eq:loss}
\end{equation}
to compute a weight $\lambda_t^{(k)}$. The total loss parameter $l_t^{(k)}$, indicates the sum of the differences between the modified desired MSE ($\sigma_m^2$) and the squared error of the first $k-1$ WLs at time $t$. Then, we add the difference $\sigma_m^2-(e_t^{(k)})^2$ to $l_t^{(k)}$, to generate $l_t^{(k+1)}$, and pass $l_t^{(k+1)}$ to the next WL, as shown in Fig. \ref{fig:boost_block}. Here, $\left[ \sigma_m^2 - \left(d_t -f_t^{(k)}(\vx_t) \right)^2 \right]$ measures how much the $k^{th}$ WL is off with respect to the final MSE performance goal. For example, in a stationary environment, if $d_t = f(\vx_t)+\nu_t$, where $f(\cdot)$ is a deterministic function and $\nu_t$ is the observation noise, one can select the desired MSE $\sigma_d^2$ as an upper bound on the variance of the noise process $\nu_t$, and define a \textit{modified} desired MSE as $\sigma_m^2 \triangleq \frac{\sigma_d^2 - \kappa}{1 - \kappa}$. In this sense,  $l_t^{(k)}$ measures how the WLs $j=1,\ldots, k$ are cumulatively performing on $(d_t,\vec{x}_t)$ pair with respect to the final performance goal.\par
We then use the weight $\lambda_t^{(k)}$ to update the $k^{th}$ WL with the ``weighted updates'', ``data reuse'', or ``random updates'' method, which we explain later in Sections \ref{sec:BRLS} and \ref{sec:BLMS}. Our aim is to make $\lambda_t^{(k)}$ large if the first $k-1$ WLs made large errors on $d_t$, so that the $k^{th}$ WL gives more importance to $(\vec{x}_t,d_t)$ in order to rectify the performance of the overall system. We now explain how to construct these weights, such that $0 < \lambda_t^{(k)} \leq 1$. To this end, we set $\lambda_t^{(1)}=1$, for all $t$, and introduce a weighting similar to (\cite{smooth_boost,onlineboost}). We define the weights as
\begin{equation}
\lambda_t^{(k)} = \min \left\{ 1, \left(\sigma^2 \right)^{l_t^{(k)}/2} \right\} \label{eq:weights_original},
\end{equation}
where $\sigma^2$ is the guaranteed upper bound on the weighted MSE of the weak learners. However, since there is no prior information about the exact MSE performance of the weak learners, we use the following weighting scheme
\begin{equation}
\lambda_t^{(k)} = \min \left\{ 1, \left(\delta_{t-1}^{(k)} \right)^{c\; l_t^{(k)}} \right\} \label{eq:weights},
\end{equation}
where $\delta_{t-1}^{(k)}$ indicates an estimate of the $k^{th}$ weak learner's MSE, and $c \geq 0$ is a design parameter, which determines the ``dependence'' of each WL update on the performance of the previous WLs, i.e., $c=0$ corresponds to ``independent'' updates, like the ordinary combination of the WLs in adaptive filtering (\cite{kozat,adamix}), while a greater $c$ indicates the greater effect of the previous WLs performance on the weight $\lambda_t^{(k)}$ of the current WL. Note that including the parameter $c$ does not change the validity of our proofs, since one can take $\left(\delta_{t-1}^{(k)} \right)^{2c}$ as the new guaranteed weighted MSE. Here, $\delta_{t-1}^{(k)}$ is an estimate of the  ``Weighted Mean Squared Error'' (WMSE) of the $k^{th}$ WL over $\{\vec{x}_t\}_{t \geq 1}$ and $\{d_t\}_{t \geq 1}$. In the basic implementation of the online boosting (\cite{smooth_boost,onlineboost}), $ \left(1-\delta_{t-1}^{(k)}\right)$ is set to the classification  advantage of the weak learners (\cite{smooth_boost}), where this advantage is assumed to be the same for all weak learners. In this paper, to avoid using any a priori knowledge and to be completely adaptive, we choose $\delta_{t-1}^{(k)}$ as the weighted and thresholded MSE of the $k^{th}$ WL up to time $t-1$ as
\begin{align}
\delta_t^{(k)} & =\frac{\ \mathlarger{\sum}_{\tau=1}^{t} \frac{\lambda_{\tau}^{(k)}}{4} \left(d_{\tau} -\left[f_{\tau}^{(k)}(\vx_{\tau})\right]^{+}\right)^2}{\sum_{\tau=1}^{t}\lambda_{\tau}^{(k)}} \nonumber \\
& = \frac{\Lambda_{t-1}^{(k)}\delta_{t-1}^{(k)}+\frac{\lambda_{t}^{(k)}}{4} \left(d_{t} -\left[f_{t}^{(k)}(\vx_{t})\right]^{+}\right)^2}{\Lambda_{t-1}^{(k)}+\lambda_{t}^{(k)}}, \label{eq:delta}
\end{align}
where $\Lambda_t^{(k)} \triangleq \sum_{\tau=1}^{t} \lambda_{\tau}^{(k)}$, and $\left[f_{\tau}^{(k)}(\vx_{\tau})\right]^{+}$ thresholds $f_{\tau}^{(k)}(\vx_{\tau})$ into the range $[-1,1]$. This thresholding is necessary to assure that $0 < \delta_t^{(k)} \leq 1$, which guarantees  $0 < \lambda_t^{(k)} \leq 1$ for all $k=1,\ldots,m$ and $t$. We point out that \eqref{eq:delta} can be recursively calculated.\par

Regarding the definition of $\lambda_t^{(k)}$, if the first $k$ WLs are ``good'', we will pass less weight to the next WLs, such that those WLs can concentrate more on the other samples. Hence, the WLs can increase the diversity by concentrating on different parts of the data \cite{kozat}. Furthermore, following this idea, in \eqref{eq:weights},  the weight $\lambda_t^{(k)}$ is larger, i.e., close to 1, if most of the WLs, $1,\ldots,k-1$, have errors larger than $\sigma_m^2$ on $(\vec{x}_t,d_t)$, and smaller, i.e., close to 0, if the pair $(\vec{x}_t,d_t)$ is easily modeled by the previous WLs such that the WLs $k,\ldots,m$ do not need to concentrate more on this pair.\par

\subsection{The Combination Algorithm \label{sec:final_rls}}
Although in the proof of our algorithm, we assume a constant combination vector $\vz$ over time, we use a time varying combination vector in practice, since there is no knowledge about the exact number of the required week learners for each problem. Hence, after $d_t$ is revealed, we also update the final combination weights $\vec{z}_t$ based on the final output $\hat{d}_t = \vec{z}_t^T \vec{y}_t$, where  $\hat{d}_t = \vz_t^T\vec{y}_t$, $\vec{y}_t = [\hat{d}_t^{(1)}, \ldots, \hat{d}_t^{(m)}]^T$. To update the final  combination weights, we use the normalized SGD algorithm \cite{sayed_book} yielding
\begin{equation}
\vz_{t+1} = \vz_{t} + \mu_z e_t \frac{\vy_t}{\|\vy_t\|^2} .
\end{equation}
%\begin{align*}
%&\vec{R}_{t+1}  = \lambda \vec{R}_{t}+ \vec{y}_t \vec{y}_t^T,  \\
%&\vec{p}_{t+1}  =  \lambda \vec{p}_{t} +  \vec{y}_t d_t,
%\end{align*}
%and
%\begin{align}
%& e_t=d_t-\vec{z}_t^T\vec{y}_{t}, \nonumber\\
%& \vec{g}_t = \frac{\vec{P}_t \vec{y}_t}{\lambda+ \vec{y}_t^T \vec{P}_{t} \vec{y}_t }, \label{eq:gdef}\\
%& \vec{z}_{t+1} = \vec{z}_t + e_t \vec{g}_t, \nonumber\\
%& \vec{P}_{t+1} = \lambda^{-1}\vec{P}_t-\lambda^{-1} \vec{g}_t \vec{y}_t^T \vec{P}_t, \nonumber
%\end{align}
%where $0<\lambda \leq 1$ is the exponential weighting.

\subsection{Choice of Parameter Values}
The choice of $\sigma_m^2$ is a crucial task, i.e., we cannot reach any desired MSE for any data sequence unconditionally. As an example, suppose that the data are generated randomly according to a known distribution, while they are contaminated with a white noise process. It is clear that we cannot obtain an MSE level below the noise power. However, if the WLs are guaranteed to satisfy the conditions of Theorem 1, this would not happen. Intuitively, there is a guaranteed upper bound (i.e., $\sigma^2$) on the worst case performance, since in the weighted MSE, the samples with a higher error have a more important effect. On the other hand, if one chooses a $\sigma_m^2$ smaller than the noise power, $l_t^{(k)}$ will be negative for almost every $k$, turning most of the weights into 1, and as a result the weak learners fail to reach a weighted MSE smaller than $\sigma^2$. Nevertheless, in practice we have to choose the parameter $\sigma_m^2$ reasonably and precisely such that the conditions of Theorem 1 are satisfied. For instance, we set $\sigma_m^2$ to be an upper bound on the noise power.\par
In addition, the number of weak learners, $m$, is chosen regarding to the computational complexity constraints. However, in our experiments we choose a moderate number of weak learners, $m=20$, which successfully improves the performance. Moreover, according to the results in Section \ref{sec:parameter}, the optimum value for $c$ is around 1, hence, we set the parameter $c=1$ in our simulations.\par
\section{Boosted NM Algorithms} \label{sec:BRLS}
At each time $t$, all of the WLs (shown in Fig. \ref{fig:boost_block}) estimate the desired data $d_t$ in parallel, and the final estimate is a linear combination of the results generated by the WLs. When the $k^{th}$ WL receives the weight $\lambda_t^{(k)}$, it updates the linear coefficients  $\vw_{t}^{(k)}$ using one of the following methods.
\subsection{Directly Using $\lambda$'s as Sample Weights \label{sec:smooth}}
Here, we consider  $\lambda_t^{(k)}$ as the weight for the observation pair $(\vec{x}_t,d_t)$ and apply a weighted NM update to $\vw_{t}^{(k)}$. For this particular weighted NM algorithm, we define the Hessian matrix and the gradient vector as
\begin{align}
&\vec{R}_{t+1}^{(k)}  \triangleq \beta \vec{R}_{t}^{(k)}+ \lambda_t^{(k)} \vec{x}_t \vec{x}_t^T, \label{eq:R} \\
&\vec{p}_{t+1}^{(k)}  \triangleq \beta \vec{p}_{t}^{(k)} + \lambda_t^{(k)} \vec{x}_t d_t, \label{eq:p}
\end{align}
where $\beta$ is the forgetting factor \cite{sayed_book} and $\vw_{t+1}^{(k)}= \left( \vec{R}_{t+1}^{(k)} \right)^{-1} \vec{p}_{t+1}^{(k)} $ can be calculated in a recursive manner as
\begin{align}
& e_t^{(k)}=d_t-\vec{x}_t^T\vw_{t}^{(k)}, \nonumber \\
& \vec{g}_{t}^{(k)} = \frac{\lambda_t^{(k)}\vec{P}_{t}^{(k)} \vec{x}_t}{\beta+\lambda_t^{(k)} \vx_t^T \vec{P}_{t}^{(k)} \vec{x}_t }, \nonumber \\
& \vw_{t+1}^{(k)} = \vw_{t}^{(k)} + e_t^{(k)} \vec{g}_{t}^{(k)}, \nonumber \\
& \vec{P}_{t+1}^{(k)} = \beta^{-1} \left(\vec{P}_{t}^{(k)}-\vec{g}_{t}^{(k)} \vec{x}_t^T \vec{P}_{t}^{(k)}\right). \label{eq:w_rls}
\end{align}
where $\vec{P}_{t}^{(k)} \triangleq \big(\vec{R}_{t}^{(k)}\big)^{-1}$, and $\vec{P}_{0}^{(k)} = v^{-1}\vI$, and $0<v \ll 1$. The complete algorithm is given in Algorithm \ref{alg:boostedWRLS} with the weighted NM implementation in \eqref{eq:w_rls}.\\ % We note that we do not need to normalize  $\lambda_t^{(k)}$ since the same weighting is used both in \eqref{eq:R} and \eqref{eq:p}. \\

\begin{algorithm}[t]
	\caption{Boosted NM-based algorithm}\label{alg:boostedWRLS}
	\begin{algorithmic}[1]
		\STATE Input: $(\vec{x}_t,d_t)$ (data stream), $m$ (number of WLs) and $\sigma_m^2$.
		\STATE Initialize the regression coefficients $\vec{w}_{1}^{(k)}$ for each WL; and the combination coefficients as $\vec{z}_1=\frac{1}{m}[1,1,\ldots,1]^T$; and for all $k$ set $\delta_0^{(k)}=0$.
		\FOR{$t=1$ \TO $T$}
		\STATE Receive the regressor data instance $\vec{x}_t$;
		%		\STATE Compute the indicator functions $s_{i,t}^{(k)}$ for all $k$'s
		\STATE Compute the WLs outputs $\hat{d}_t^{(k)}= \vec{x}_t^T\vec{w}_{t}^{(k)}$;
		\STATE Produce the final estimate $\hat{d}_t = \vec{z}_t^T[\hat{d}_t^{(1)}, \ldots, \hat{d}_t^{(m)}]^T$;
		\STATE Receive the true output $d_t$ (desired data);
		\STATE $\lambda_t^{(1)}=1$; $l_t^{(1)}=0$;
		\FOR{$k=1$ \TO $m$}
		\STATE $\lambda_t^{(k)} = \min \left\{ 1, \left(\delta_{t-1}^{(k)} \right)^{c \; l_t^{(k)}} \right\}$;
		\STATE Update the regression coefficients $\vec{w}_{t}^{(k)}$ by using the NM and the weight $\lambda_t^{(k)}$ based on one of the introduced algorithms in Section \ref{sec:BRLS};
		\STATE $e_t^{(k)}= d_t -\hat{d}_t^{(k)}$; %(the error after the update).
		\STATE $\delta_t^{(k)} = \frac{\Lambda_{t-1}^{(k)}\delta_{t-1}^{(k)}+\frac{\lambda_{t}^{(k)}}{4} \left(d_{t} -\left[f_{t}^{(k)}(\vx_{t})\right]^{+}\right)^2}{\Lambda_{t-1}^{(k)}+\lambda_{t}^{(k)}}$;
		\STATE $\Lambda_{t}^{(k)} = \Lambda_{t-1}^{(k)}+\lambda_{t}^{(k)}$
		\STATE $l_t^{(k+1)} =  l_t^{(k)}+\left[ \sigma_m^2 - \left(e_t^{(k)}\right)^2 \right]$;
		\ENDFOR
		\STATE $e_t=d_t-\vec{z}_t^T\vec{y}_{t}$;
		%		\STATE $\vec{g}_t = \frac{\vec{P}_t \vec{y}_t}{\lambda+ \vec{y}_t^T \vec{P}_{t} \vec{y}_t }$;
		%		\STATE $\vec{z}_{t+1} = \vec{z}_t + e_t \vec{g}_t$;
		%		\STATE $\vec{P}_{t+1} = \lambda^{-1}\vec{P}_t-\lambda^{-1} \vec{g}_t \vec{y}_t^T \vec{P}_t$;
		\STATE $\vec{z}_{t+1}  = \vec{z}_t+\mu_z e_t \frac{\vy_t}{\|\vy_t\|^2}$; 
		%		(update for the combination weights), where $\vec{y}_t = [\hat{d}_t^{(1)}, \ldots, \hat{d}_t^{(m)}]^T$.
		
		\ENDFOR
	\end{algorithmic}
	\label{alg:RLS}
\end{algorithm}

\noindent
%{\bf Remark 1:} We emphasize that one can be inclined to use a single boosted NM algorithm instead of running $m$ NM algorithms in parallel due to the computational complexity considerations. In this single NM implementation, the calculated importance weight can be used as a weight in the next time instant $t+1$ so that a single NM algorithm may concentrate more on wrongly regressed samples. However, running $m$ NM filters with boosting provides ``diversity'' \cite{kozat}, where each constituent algorithm concentrates on those parts of the data that other WLs are not successful, due to the boosting weights in \eqref{eq:weights}. Hence, by this boosting and then final mixture-of-experts combination \cite{kozat}, we achieve a significantly improved learning performance.

\subsection{Data Reuse Approaches Based on The Weights}

Another approach follows Ozaboost (\cite{ozaboost}). In this approach, from  $\lambda_t^{(k)}$, we generate an integer, say  $n_t^{(k)} = \mathrm{ceil}(K\lambda_t^{(k)})$, where $K$ is a design parameter that takes on positive integer values. We then apply the NM update on the $(\vec{x}_t,d_t)$ pair repeatedly $n_t^{(k)}$ times, i.e., run the NM update on the same $(\vec{x}_t,d_t)$ pair  $n_t^{(k)}$ times consecutively. Note that $K$ should be determined according to the computational complexity constraints. However, increasing $K$ does not necessarily result in a better performance, therefore, we use moderate values for $K$, e.g., we use $K=5$ in our simulations. The final $\vw_{t+1}^{(k)}$ is calculated after  $n_{t}^{(k)}$ NM updates. As a major advantage, clearly, this reusing approach can be readily generalized to other adaptive algorithms in a straightforward manner.\par
%\noindent
%{\bf Remark:} We emphasize that the data reuse approach is exactly equal to the approach in Section~\ref{sec:smooth}, where the cross correlation vector and auto correlation matrix are defined as
%\begin{align*}
%  &\vec{R}_{t+1}^{(k)}  = \vec{R}_{t}^{(k)}+ n_t^{(k)} \vec{x}_t \vec{x}_t^T,  \\
%  &\vec{p}_{t+1}^{(k)}  =  \vec{p}_{t}^{(k)} + n_t^{(k)} \vec{x}_t d_t.
%\end{align*}
%Hence, for the RLS recursion, these two approaches are equivalent up to $\mathrm{ceil}()$ operation to generate an integer from $\lambda_t^{(k+1)}$. However, such an equivalence is not usually correct for other type of updates such as the LMS algorithm, hence we provide the LMS algorithm for completeness.  \\
We point out that Ozaboost (\cite{ozaboost}) uses a different data reuse strategy. In this approach, $\lambda_t^{(k)} $ is used as the parameter of a Poisson distribution and an integer $n_{t}^{(k)}$ is randomly generated from this  Poisson distribution. One then applies the NM update $n_{t}^{(k)}$ times.\par 

\subsection{Random Updates Approach Based on The Weights}\label{sec:rndRLS}
In this approach, we simply use the weight $\lambda_t^{(k)}$ as a probability of updating the $k^{th}$ WL at time $t$. To this end, we generate a Bernoulli random variable, which is $1$ with probability $\lambda_t^{(k)}$ and is $0$ with probability $1-\lambda_t^{(k)}$. Then, we update the $k^{th}$ WL, only if the Bernoulli random variable equals $1$. With this method, we significantly reduce the computational complexity of the algorithm. Moreover, due to the dependence of this Bernoulli random variable on the performance of the previous constituent WLs, this method does not degrade the MSE performance, while offering a considerably lower complexity, i.e., when the MSE is low, there is no need for further updates, hence, the probability of an update is low, while this probability is larger when the MSE is high.\par

\section{Boosted SGD Algorithms} \label{sec:BLMS}
In this case, as shown in Fig. \ref{fig:boost_block}, we have $m$ parallel running WLs, each of which is updated using the SGD algorithm. Based on the weights given in \eqref{eq:weights} and the total loss and MSE parameters in \eqref{eq:loss} and \eqref{eq:delta}, we next introduce three SGD based boosting algorithms, similar to those introduced in Section \ref{sec:BRLS}.

\subsection{Directly Using $\lambda$'s to Scale The Learning Rates \label{sec:smooth_lms}}
We note that by construction method in \eqref{eq:weights}, $0<\lambda_t^{(k)} \leq 1$, thus, these weights can be directly used to scale the learning rates for the SGD updates. When the $k^{th}$ WL receives the weight $\lambda_t^{(k)}$, it updates its coefficients  $\vw_{t}^{(k)}$, as
\begin{equation}
\vw_{t+1}^{(k)}  = \left(\vec{I}-\mu^{(k)} \lambda_t^{(k)} \vec{x}_t \vec{x}_t^T\right) \vw_{t}^{(k)}+\mu^{(k)}\lambda_t^{(k)} \vec{x}_t d_t, \label{eq:blms}
\end{equation}
where $0<\mu^{(k)} \lambda_t^{(k)} \leq \mu^{(k)}$. Note that we can choose $\mu^{(k)} = \mu$ for all $k$, since the online algorithms work consecutively from top to bottom, and the $k^{th}$ WL will have a different learning rate $\mu^{(k)} \lambda_t^{(k)}$.\par

\subsection{A Data Reuse Approach Based on The Weights}

In this scenario, for updating $\vw_{t}^{(k)}$, we  use the SGD update $n_t^{(k)} = \mathrm{ceil}(K\lambda_t^{(k)})$ times to obtain the $\vw_{t+1}^{(k)}$ as
\begin{align}
&\vq^{(0)}  = \vw_{t}^{(k)}, \label{eq:iterate} \nonumber \\
&\vq^{(a)}  = \left(\vec{I}-\mu^{(k)} \vec{x}_t \vec{x}_t^T\right) \vq^{(a-1)}+\mu^{(k)} \vec{x}_t d_t, \; a=1,\ldots,n_{t}^{(k)},  \nonumber \\
& \vw_{t+1}^{(k)} = \vq^{\left(n_{t}^{(k)}\right)}. 
\end{align}
where $K$ is a constant design parameter.\par
Similar to the NM case, if we follow the Ozaboost (\cite{ozaboost}), we use the weights to generate a random number $n_t^{(k)}$ from a Poisson distribution with parameter $\lambda_t^{(k)}$, and perform the SGD update $n_t^{(k)}$ times on $\vw_{t}^{(k)}$ as explained above.\par

\subsection{Random Updates Based on The Weights} \label{sec:rndLMS}
Again, in this scenario, similar to the NM case, we use the weight $\lambda_t^{(k)}$ to generate a random number from a Bernoulli distribution, which equals $1$ with probability $\lambda_t^{(k)}$, and equals $0$ with probability $1-\lambda_t^{(k)}$. Then we update $\vw_t$ using SGD only if the generated number is $1$.

\section{Analysis Of The Proposed Algorithms}\label{sec:analysis}
In this section we provide the complexity analysis for the proposed algorithms. We prove an upper bound for the weights $\lambda_t^{(k)}$, which is significantly less than 1. This bound shows that the complexity of the ``random updates'' algorithm is significantly less than the other proposed algorithms, and slightly greater than that of a single WL. Hence, it shows the considerable advantage of ``boosting with random updates'' in processing of high dimensional data. 

%Furthermore, we use this bound in MSE analysis of the algorithms. In addition, for the sake of simplicity, we have chosen the ``dependence parameter'' $c = 1$. Nevertheless, the results can be easily extended to the general case.
%\vspace{-0.4cm}
\subsection{Complexity Analysis}
Here we compare the complexity of the proposed algorithms and find an upper bound for the computational complexity of random updates scenario (introduced in Section \ref{sec:rndRLS} for NM, and in Section \ref{sec:rndLMS} for SGD updates), which shows its significantly lower computational burden with respect to two other approaches. For $\vec{x}_t \in \Real^r$, each WL performs $O(r)$ computations to generates its estimate, and if updated using the NM algorithm, requires $O(r^2)$ computations due to updating the matrix $\bR_{t}^{(k)}$, while it needs $O(r)$ computations when updated using the SGD method (in their most basic implementation).\par
We first derive the computational complexity of using the NM updates in different boosting scenarios. Since there are a total of $m$ WLs, all of which are updated in the ``weighted updates'' method, this method has a computational cost of order $O(mr^2)$ per each iteration $t$. However, in the ``random updates'', at iteration $t$, the $k^{th}$ WL may or may not be updated with probabilities $\lambda_t^{(k)}$ and $1-\lambda_t^{(k)}$ respectively, yielding
\begin{equation}
C_t^{(k)} = 
\begin{cases}
O(r^2) & \quad \text{with probability}\  \lambda_t^{(k)}\\
O(r) & \quad \text{with probability}\  1-\lambda_t^{(k)},
\end{cases}
\end{equation}
where $C_t^{(k)}$ indicates the complexity of running the $k^{th}$ WL at iteration $t$. Therefore, the total computational complexity $C_t$ at iteration $t$ will be $C_t = \sum_{k=1}^{m} C_t^{(k)}$, which yields
\begin{equation}
E\left[C_t\right] = E\left[\sum_{k=1}^{m} C_t^{(k)}\right] = \sum_{k=1}^{m} E[\lambda_t^{(k)}] O(r^2)
\end{equation}
Hence, if $E\big[\lambda_t^{(k)}\big]$ is upper bounded by $\tilde{\lambda}^{(k)} < 1$, the average computational complexity of the random updates method, will be
\begin{equation}
E\left[C_t\right] < \sum_{k=1}^{m} \tilde{\lambda}^{(k)} O(r^2).
\end{equation}
In Theorem 2, we provide sufficient constraints to have such an upper bound.\par
Furthermore, we can use such a bound for the ``data reuse'' mode as well. In this case, for each WL $f_t^{(k)}$, we perform the NM update $\lambda_t^{(k)}K$ times, resulting a computational complexity of order $\displaystyle E\left[C_t\right] < \sum_{k=1}^{m} K\ \tilde{\lambda}^{(k)} (O(r^2))$. For the SGD updates, we similarly obtain the computational complexities $O(mr)$, $\sum_{k=1}^{m} O\big(\tilde{\lambda}^{(k)}r\big)$, and $ \sum_{k=1}^{m} O\big(K\tilde{\lambda}^{(k)}r\big)$, for the ``weighted updates'', ``random updates'', and ``data reuse'' scenarios respectively.\par
The following theorem determines the upper bound $\tilde{\lambda}^{(k)}$ for $E\big[\lambda_t^{(k)}\big]$.\\
{\bf Theorem 2.} {\em If the WLs converge and achieve a sufficiently small MSE (according to the proof following this Theorem), the following upper bound is obtained for $\lambda_t^{(k)}$, given that $\sigma_m^2$ is chosen properly,
	\begin{equation}\label{eq:bound}
	E \left[\lambda_t^{(k)}\right] \leq \tilde{\lambda}^{(k)} = \left(\gamma^{-2 \sigma_m^2}(1+2 \zeta^2 \ln \gamma)\right)^{\frac{1-k}{2}},\\
	\end{equation}
	where $\gamma \triangleq E\left[\delta_{t-1}^{(k)}\right]$ and $\zeta^2 \triangleq E\left[\left(e_t^{(k)}\right)^2\right]$.}\\
It can be straightforwardly shown that, this bound is less than $1$ for appropriate choices of $\sigma_m^2$, and reasonable values for the MSE according to the proof. This theorem states that if we adjust $\sigma_m^2$ such that it is achievable, i.e., the WLs can provide a slightly lower MSE than $\sigma_m^2$, the probability of updating the WLs in the random updates scenario will decrease. This is of course our desired results, since if the WLs are performing sufficiently well, there is no need for additional updates. Moreover, if $\sigma_m^2$ is opted such that the WLs cannot achieve a MSE equal to $\sigma_m^2$, the WLs have to be updated at each iteration, which increases the complexity.\par
\noindent
{\bf Proof:} For simplicity, in this proof, we have assumed that $c=1$, however, the results are readily extended to the general values of $c$. We construct our proof based on the following assumption:\\
{\bf Assumption:} assume that $e_t^{(k)}$'s are independent and identically distributed (i.i.d) zero-mean Gaussian random variables with variance $\zeta^2$.\\
We have
\begin{align}
E\left[\lambda_t^{(k)}\right] & = E\left[\min \left\{1,\left(\delta_{t-1}^{(k)}\right)^{l_t^{(k)}} \right\}\right]\nonumber\\
& \leq \min\left\{1,E\left[\left(\delta_{t-1}^{(k)}\right)^{l_t^{(k)}}\right]\right\}
\end{align}
Now, we show that under certain conditions, $E\big[\big(\dltk \big)^{\ltk} \big]$ will be less than 1, hence, we obtain an upper bound for $E\big[\lmtk\big]$. We define $s \triangleq \ln(\dltk)$, yielding
\begin{align}
E\left[\left(\dltk \right)^{\ltk} \right] = E \left[E\left[ \exp\big(s\ {\ltk}\big) \Big|s \right]\right] = E\left[M_{\ltk}(s) \Big| s \right],
\end{align}
where $M_{\ltk}(.)$ is the moment generating function of the random variable $\ltk$.
From the Algorithm \ref{alg:RLS}, $\ltk = (k-1)\sigma_m^2-\sum_{j=1}^{k-1}\big(e_t^{(j)} \big)^2$.
According to the Assumption, $\frac{e_t^{(j)}}{\zeta}$ is a standard normal random variable. Therefore, $\sum_{j=1}^{k-1}\big(e_t^{(j)} \big)^2$ has a Gamma distribution as $\Gamma\big(\frac{k-1}{2},2\zeta^2\big)$ (\cite{papoulis_book}), which results in the following moment generating function for $\ltk$
\begin{align}
M_{\ltk}(s) & = \exp\left(s(k-1)\sigma_m^2\right)\left(1+2\zeta^2s\right)^{\frac{1-k}{2}} \nonumber \\
& = \left(\dltk\right)^{(k-1)\sigma_m^2} \left(1+2\zeta^2 \ln \left(\dltk\right)\right)^{\frac{1-k}{2}}. \label{eq:moment}
\end{align}
In the above equality $\dltk$ is a random variable, the mean of which is denoted by $\gamma$. We point out that $\gamma$ will approach to $\zeta^2$ in convergence. We define a function $\varphi(.)$ such that $E\left[\lmtk\right] = E\left[\varphi\left(\dltk\right)\right]$, and seek to find a condition for $\varphi(.)$ to be a concave function. Then, by using the Jenssen's inequality for concave functions, we have
\begin{equation}\label{eq:Jen}
E\left[\lmtk\right] \leq \varphi(\gamma).
\end{equation}
Inspired by \eqref{eq:moment}, we define $A\left(\dltk\right) \triangleq {\dltk}^{-2\sigma_m^2}\left(1+2\zeta^2\ln\left(\dltk\right)\right)$ and $\varphi\left(\dltk\right) \triangleq \left(A\left(\dltk\right)\right)^{\frac{1-k}{2}}$. By these definitions we obtain
\begin{small}
	\begin{multline}
	\hspace{-0.5cm}
	\varphi''\left(\dltk\right) = \frac{1-k}{2}\left(A\left(\dltk\right)\right)^{\frac{-k-3}{2}}\Bigg[\left(\frac{-k-1}{2}\right)\left(A'\left(\dltk\right)\right)^2 \\ +\left(A\left(\dltk\right)\right)^2 A''\left(\dltk\right)\Bigg].
	\end{multline}
\end{small}
Considering that $k>1$, in order for $\varphi(.)$ to be concave, it suffices to have
\begin{equation}
\left(A\left(\dltk\right)\right)^2 A''\left(\dltk\right) > \left(\frac{k+1}{2}\right)\left(A'\left(\dltk\right)\right)^2,
\end{equation}
which reduces to the following necessary and sufficient conditions:
\begin{equation}
\frac{\left(\dltk\right)^{2\sigma_m^2}}{\left(1+2\zeta^2\ln\left(\dltk\right)\right)^2} < \frac{\left(1+2\sigma_m^2\right)^2}{4(k+1)},
\end{equation}
and
\begin{equation}
\frac{(1-\xi_1)\sigma_m^2}{1-2\sigma_m^2\ln\left(\dltk\right)}<\zeta^2<\frac{(1-\xi_2)\sigma_m^2}{1-2\sigma_m^2\ln\left(\dltk\right)},
\end{equation}
where
\[
\xi_1 = \frac{\alpha^2(1+2\sigma_m^2)+\alpha\sqrt{(1+2\sigma_m^2)^2\alpha^2-4(k+1)(\dltk)^{2\sigma_m^2}}}{2(k+1)(\dltk)^{2\sigma_m^2}},\\
\]
\[
\xi_2 = \frac{\alpha^2(1+2\sigma_m^2)-\alpha\sqrt{(1+2\sigma_m^2)^2\alpha^2-4(k+1)(\dltk)^{2\sigma_m^2}}}{2(k+1)(\dltk)^{2\sigma_m^2}},
\]
and
\[
\alpha \triangleq 1+2\zeta^2\ln\left(\dltk\right).
\]
Under these conditions, $\varphi(.)$ is concave, therefore, by substituting $\varphi(.)$ in \eqref{eq:Jen} we achieve \eqref{eq:bound}. This concludes the proof of the Theorem 2. $\Box$\par

\section{Experiments}
In this section, we demonstrate the efficacy of the proposed boosting algorithms for NM and SGD linear WLs under different scenarios. To this end, we first consider the ``online regression'' of data generated with a stationary linear model. Then, we illustrate the performance of our algorithms under nonstationary conditions, to thoroughly test the adaptation capabilities of the proposed boosting framework. Furthermore, since the most important parameters in the proposed methods are $\sigma_m^2$, $c$, and $m$, we investigate their effects on the final MSE performance. Finally, we provide the results of the experiments over several real and synthetic benchmark datasets. \par
Throughout this section, ``SGD'' represents the linear SGD-based WL, ``NM'' represents the linear NM-based WL, and a prefix ``B'' indicates the boosting algorithms. In addition, we use the suffixes ``-WU'', ``-RU'', or ``-DR'' to denote the ``weighted updates'', ``random updates'', or ``data reuse'' modes, respectively, e.g., the ``BSGD-RU'' represents the ``Boosted SGD-based algorithm using Random Updates''.\par
In order to observe the boosting effect, in all experiments, we set the step size of SGD and the forgetting factor of the NM to their optimal values, and use those parameters for the WLs, too. In addition, the initial values of all of the weak learners in all of the experiments are set to zero. However, in all experiments, since we use $K=5$ in BSGD-DR algorithm, we set the step size of the WLs in BSGD-DR method to $\mu / K= \mu/5$, where, $\mu$ is the step size of the SGD. To compare the MSE results, we have provided the Accumulated Square Error (ASE) results.

\subsection{Stationary Data}
In this experiment, we consider the case where the desired data is generated by a stationary linear model. The input vectors $\vx_t = [x_1\ x_2\ 1]$ are 3-dimensional, where $[x_1\ x_2]$ is drawn from a jointly Gaussian random process and then scaled such that $\vx_t=[x_1\ x_2]^T \in [0\ 1]^2$. We include 1 as the third entry of $\vx_t$ to consider affine learners. Specifically the desired data is generated by $d_t = [1\ 1\ 1]^T\ \vx_t + \nu_t$, where $\nu_t$ is a random Gaussian noise with a variance of $0.01$.\par
In our simulations, we use $m=20$ WLs and $\mu=0.1$ for all SGD learners. In addition, for NM-based boosting algorithms, we set the forgetting factor $\beta = 0.9999$ for all algorithms. Moreover, we choose $\sigma_m^2 = 0.02$ for SGD-based algorithms and $\sigma_m^2 = 0.004$ for NM-based algorithms, $K=5$ for data reuse approaches, and $c=1$ for all boosting algorithms. To achieve robustness, we average the results over 100 trials.\par
\begin{figure}
	\centering
	\includegraphics[width=0.6\textwidth]{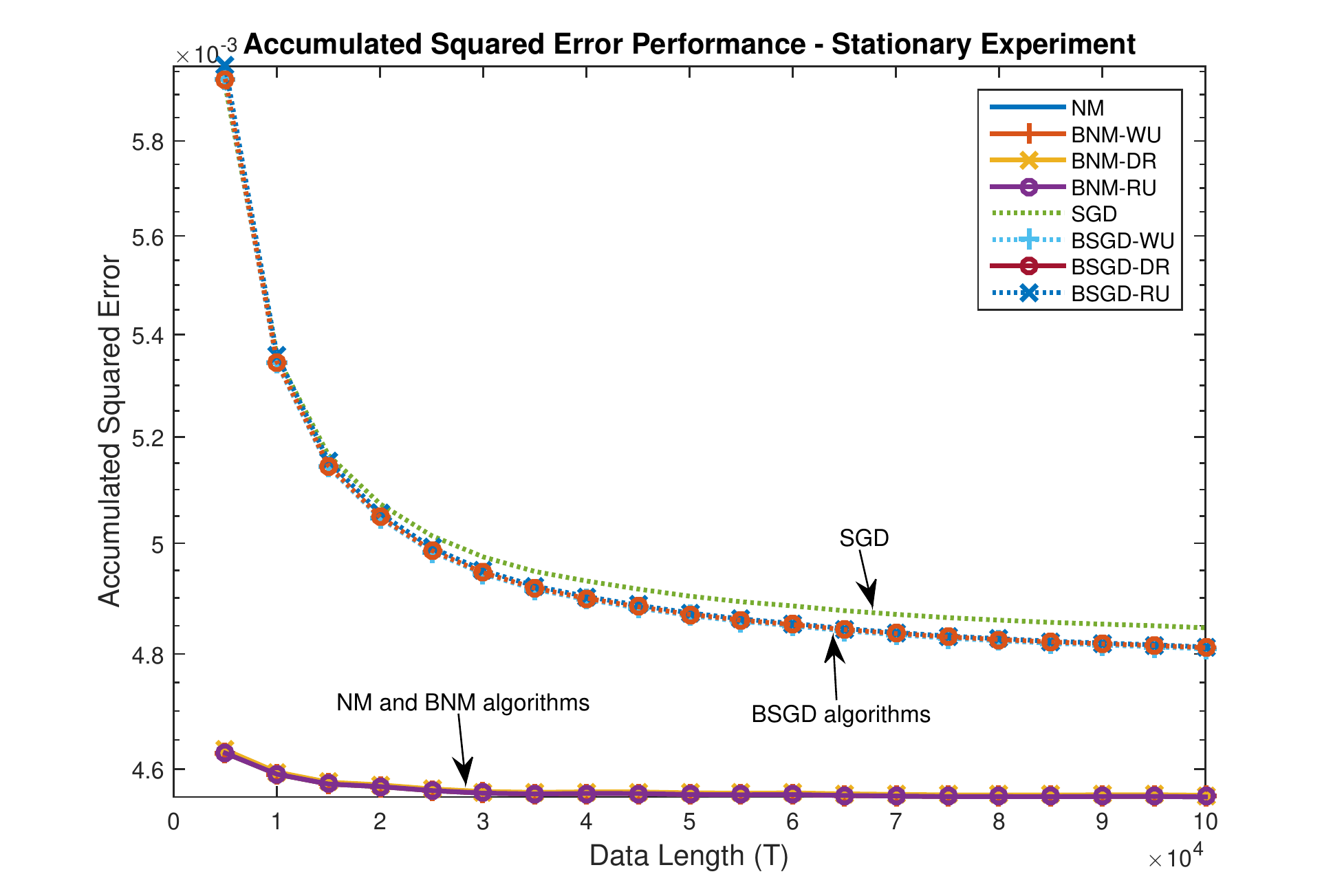}\\
	\caption{The ASE performnce of the proposed algorithms in the stationary data experiment.} \label{fig:stationary}
\end{figure}

As depicted in Fig. \ref{fig:stationary}, our proposed methods boost the performance of a single linear SGD-based WL. Nevertheless, we cannot further improve the performance of a linear NM-based WL in such a stationary experiment since the NM achieves the lowest MSE. We point out that the random updates method achieves the performance of the weighted updates method and the data reuse method with a much lower complexity. In addition, we observe that by increasing the data length, the performance improvement increases (Note that the distance between the ASE curves is slightly increasing).\par

%\begin{figure}
%	\centering
%	\includegraphics[width=0.4\textwidth]{figures/new_results/sta_LMS_relimp.eps}\\
%	\caption{The relative improvement in ASE performnce of the RLS-based algorithms in the stationary data experiment.} \label{fig:sta_lms_rel}
%\end{figure}

\subsection{Chaotic Data} \label{sec:NonSt}
\begin{figure}
	\centering
	\includegraphics[width=0.6\textwidth]{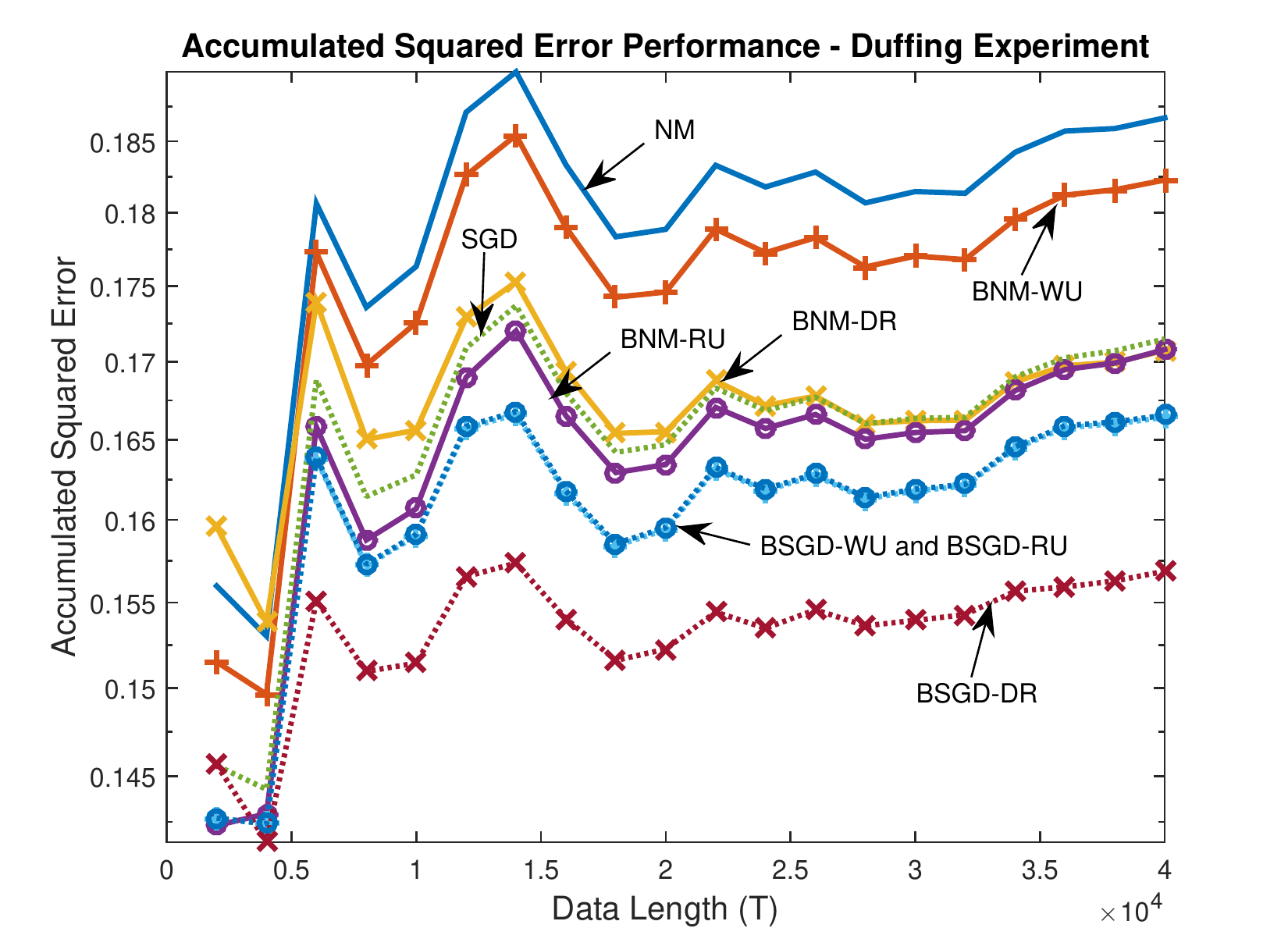}\\
	\caption{ASE performance of the proposed methods on a Duffing data set.}
	\label{fig:duffing}
	
\end{figure}

Here, in order to show the tracking capability of our algorithms in nonstationary environments, we consider the case where the desired data is generated by the Duffing map (\cite{wiggins2003}) as a chaotic model. Specifically, the data is generated by the following equation $x_{t+1} = 2.75 x_t - x_t^3 - 0.2 x_{t-1}$, where we set $x_{-1} = 0.9279$ and $x_0 = 0.1727$. We consider $d_t = x_{t+1}$ as the desired data and $[x_{t-1}\ x_t\ 1]$ as the input vector.
%for $1 \leq t < T/4$
%\begin{equation*}
%d_t = \begin{cases}
%[1\ 1]\vx_t+\nu_t &\quad \text{if}\ \Phi_0^T \vx_t < 0.8\\
%0.5+\nu_t &\quad \text{if}\ 0.8 \leq \Phi_0^T \vx_t < 1.2\\
%[1\ -0.5]\vx_t+\nu_t &\quad \text{if}\ \Phi_0^T \vx_t \geq 1.2\\
%\end{cases}
%\end{equation*}
%for $T/4 \leq t < T/2$
%\begin{equation*}
%d_t = \begin{cases}
%[1\ -0.5]\vx_t+\nu_t &\quad \text{if}\ \Phi_0^T \vx_t < 0.7\\
%[1\ 1]\vx_t+\nu_t &\quad \text{if}\ 0.7 \leq \Phi_0^T \vx_t < 1.1\\
%[-1\ 0.5]\vx_t+\nu_t &\quad \text{if}\ \Phi_0^T \vx_t \geq 1.1\\
%\end{cases}
%\end{equation*}
%for $T/2 \leq t < 3T/4$
%\begin{equation*}
%d_t = \begin{cases}
%[0\ 1.5]\vx_t+\nu_t &\quad \text{if}\ \Phi_0^T \vx_t < 0.5\\
%[1\ -0.5]\vx_t+\nu_t &\quad \text{if}\ 0.5 \leq \Phi_0^T \vx_t < 1.5\\
%[2\ 0.5]\vx_t+\nu_t &\quad \text{if}\ \Phi_0^T \vx_t \geq 1.5\\
%\end{cases}
%\end{equation*}
%for $3T/4 \leq t \leq T$
%\begin{equation*}
%d_t = \begin{cases}
%[4\ -0.5]\vx_t+\nu_t &\quad \text{if}\ \Phi_0^T \vx_t < 0.8\\
%[2\ 0.5]\vx_t+0.8+\nu_t &\quad \text{if}\ 0.8 \leq \Phi_0^T \vx_t < 1.2\\
%[1\ 1]\vx_t+\nu_t &\quad \text{if}\ \Phi_0^T \vx_t \geq 1.2\\
%\end{cases}
%\end{equation*}
In this experiment, each boosting algorithm uses 20 WLs. The step sizes for the SGD-based  algorithms are set to $0.1$, the forgetting factor $\beta$ for the NM-based algorithms are set to $0.999$, and the modified desired MSE parameter $\sigma_m^2$ is set to $0.25$ for BSGD methods, and $0.17$ for the BNM methods. Note that although the value of $\sigma_m^2$ is higher than the achieved MSE, it can improve the performance significantly. This is because of the boosting effect, i.e., emphasizing on the harder data patterns. The figures show the superior performance of our algorithms over a single WL (whose step size is chosen to be the best), in this highly nonstationary environment. Moreover, as shown in Fig. \ref{fig:duffing}, in the SGD-based boosted algorithms, the data reuse method shows a better performance relative to the other boosting methods. However, the random updates method has a significantly lower time consumption, which makes it desirable for larger data lengths. From the Fig. \ref{fig:duffing}, one can see that our method is truly boosting the performance of the conventional linear WLs in this chaotic environment.\par
\begin{figure}
	\centering
	\includegraphics[width=0.6\textwidth]{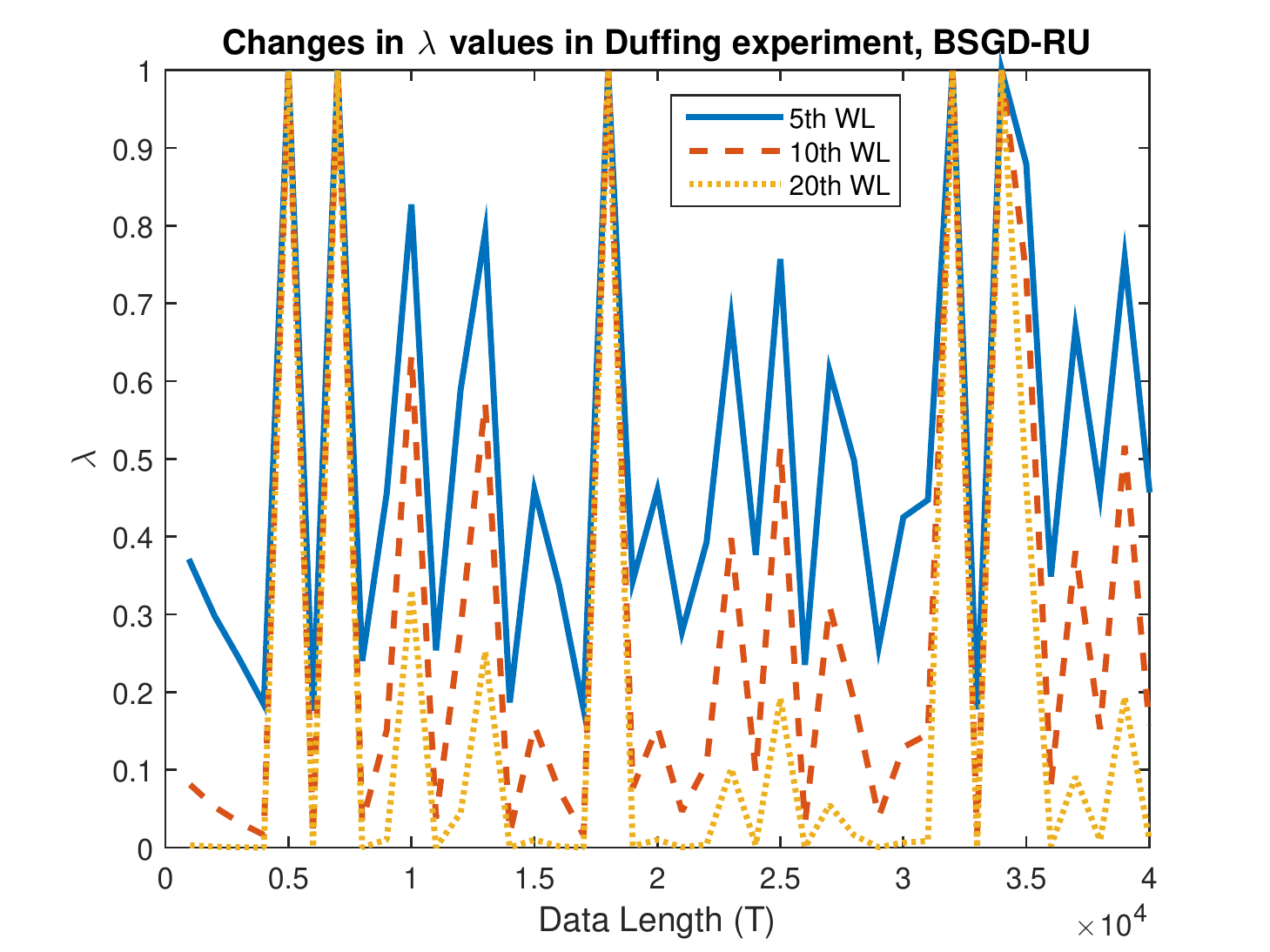}\\
	\caption{The changing of the weights in BSGD-RU algorithm in the Duffing data experiment.}\label{fig:landa}
\end{figure}

From the Fig. \ref{fig:landa}, we observe the approximate changes of the weights, in the BSGD-RU algorithm running over the Duffing data. As shown in this figure, the weights do not change monotonically, and this shows the capability of our algorithm in effective tracking of the nonstationary data. Furthermore, since we update the WLs in an ordered manner, i.e., we update the $(k+1)^{th}$ WL after the $k^{th}$ WL is updated, the weights assigned to the last WLs are generally smaller than the weights assigned to the previous WLs. As an example, in Fig. \ref{fig:landa} we see that the weights assigned to the $5^{th}$ WL are larger than those of the $10^{th}$ and $20^{th}$ WLs. Furthermore, note that in this experiment, the dependency parameter $c$ is set to 1. We should mention that increasing the value of this parameter, in general, causes the lower weights, hence, it can considerably reduce the complexity of the random updates and data reuse methods.

\subsection{The Effect of Parameters}\label{sec:parameter}
\begin{figure*}
	\centering
	\begin{subfigure}[b]{0.5\textwidth}
		\centering
		\includegraphics[scale=0.6]{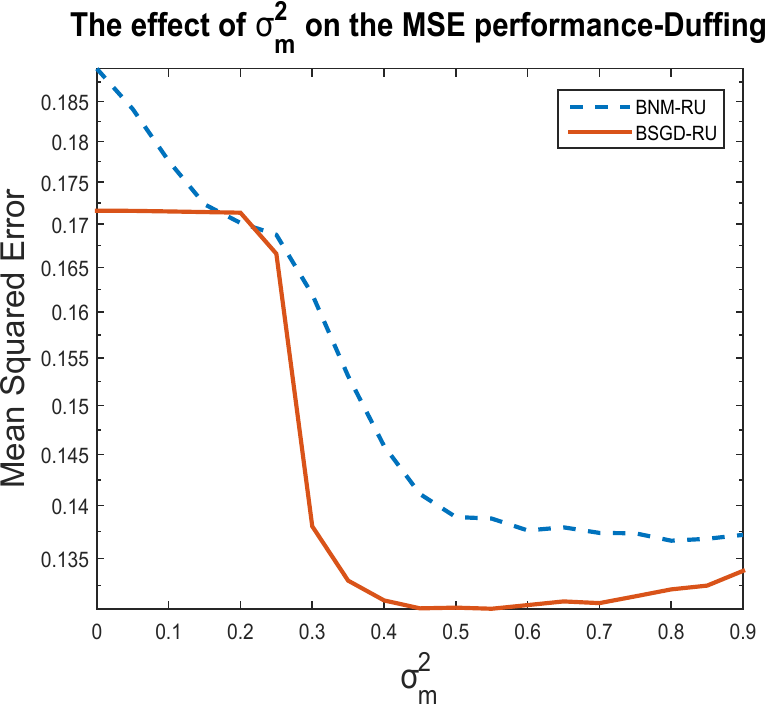}
		\caption{The effect of the parameter $\sigma_m^2$}\label{fig:sigma}
	\end{subfigure}
	\begin{subfigure}[b]{0.5\textwidth}
		\centering
		\includegraphics[scale=0.6]{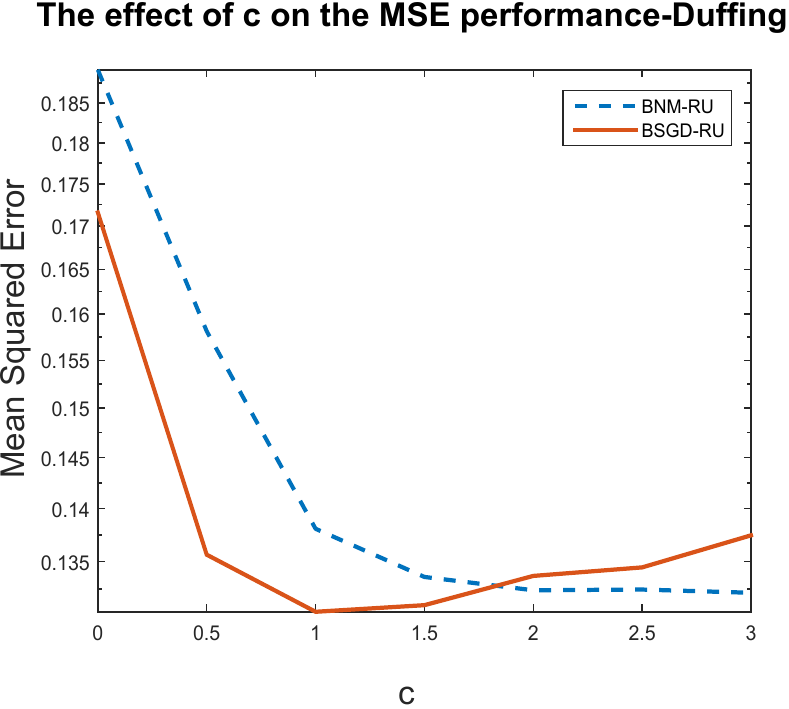}
		\caption{The effect of the parameter $c$}\label{fig:c}
	\end{subfigure}
	\begin{subfigure}[b]{0.5\textwidth}
		\centering
		\includegraphics[scale=0.6]{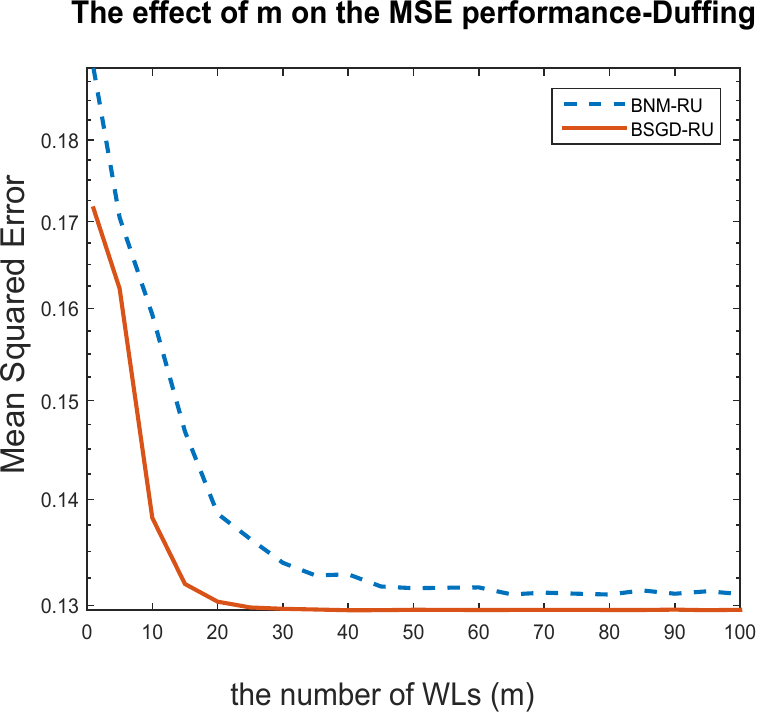}
		\caption{The effect of the parameter $m$}\label{fig:m}
	\end{subfigure}
	%	\begin{subfigure}[b]{0.4\textwidth}
	%		\centering
	%		\includegraphics[scale=0.5]{figures/new_results/c_LMS_relimp.eps}
	%		\caption{The effect of the parameter $c$}\label{fig:c_lms}
	%	\end{subfigure}
	\caption{The effect of the parameters $\sigma_m^2$, $c$, and $m$, on the MSE performance of the BNM-RU and BSGD-RU algorithms in the Duffing data experiment.}\label{fig:parameter}
\end{figure*}

In this section, we investigate the effects of the dependence parameter $c$ and the modified desired MSE $\sigma_m^2$ as well as the number of WLs ,$m$, on the boosting performance of our methods in the Duffing data experiment, explained in Section \ref{sec:NonSt}. From the results in Fig. \ref{fig:m}, we observe that, increasing the number of WLs up to $30$ can improve the performance significantly, while further increasing of $m$ only increases the computational complexity without improving the performance. In addition, as shown in Fig. \ref{fig:c}, in this experiment, the dependency parameter $c$ has an optimum value around $1$. We note that choosing small values for $c$ reduces the boosting effect, and causes the weights to be larger, which in turn increases the computational complexity in random updates and data reuse approaches. On the other hand, choosing very large values for $c$ increases the dependency, i.e., in this case the generated weights are very close to $1$ or $0$, hence, the boosting effect is decreased. Overall, one should choose values around $1$ for $c$ to avoid those extreme cases.\par
Furthermore, as depicted in Fig. \ref{fig:sigma}, there is an optimum value around $0.5$ for $\sigma_m^2$ in this experiment. Note that, choosing small values for $\sigma_m^2$ results in large weights, thus, increases the complexity and reduces the diversity. However, choosing higher values for $\sigma_m^2$ results in smaller weights, and in turn reduces the complexity. Nevertheless, we note that increasing the value of $\sigma_m^2$ does not necessarily enhance the performance. Through the experiments, we find out that $\sigma_m^2$ must be in the order of the MSE amount to obtain the best performance.

%\begin{figure}
%	\centering
%	\includegraphics[width=0.3\textwidth]{figures/new_results/m_RLS_relimp.eps}\\
%	\caption{The effect of the parameter $m$ on the relative performance improvement of the RLS-based algorithms in the nonstationary data experiment.}\label{fig:m_rls}
%\end{figure}
%
%\begin{figure}
%	\centering
%	\includegraphics[width=0.3\textwidth]{figures/new_results/m_LMS_relimp.eps}\\
%	\caption{The effect of the parameter $c$ on the relative performance improvement of the LMS-based algorithms in the nonstationary data experiment.}\label{fig:m_lms}
%\end{figure}
%
%\begin{figure}
%	\centering
%	\includegraphics[width=0.3\textwidth]{figures/new_results/c_RLS_relimp.eps}\\
%	\caption{The effect of the parameter $c$ on the relative performance improvement of the RLS-based algorithms in the nonstationary data experiment.}\label{fig:c_rls}
%\end{figure}
%
%\begin{figure}
%	\centering
%	\includegraphics[width=0.3\textwidth]{figures/new_results/c_LMS_relimp.eps}\\
%	\caption{The effect of the parameter $c$ on the relative performance improvement of the LMS-based algorithms in the nonstationary data experiment.}\label{fig:c_lms}
%\end{figure}

\subsection{Benchmark Real and Synthetic Data Sets}
In this section, we demonstrate the efficiency of the introduced methods over some widely used real life machine learning regression data sets. We have normalized each dimension of the data to the interval $[-1,1]$ in all algorithms. We present the MSE performance of the algorithms in Table \ref{tab:result}. These experiments show that our algorithms can successfully improve the performance of single linear WLs. We now describe the experiments and provide the results:\\

\begin{table*}[t!]
	\tiny
	\centering
	\begin{tabular} {|c|c|c|c|c||c|c|c|c|} \hline
		\backslashbox[3em]{\bf{Data Sets} \kern-2em}{\kern-2em \bf{Algorithms}}
		& SGD & BSGD-WU & BSGD-DR & BSGD-RU & NM & BNM-WU & BNM-DR & BNM-RU
		\\ \hline\hline
		MV & 0.2711 & 0.2707 & 0.2706 & 0.2707 & 0.2592 & 0.2645 & 0.2587 & 0.2584  \\ \hline
		\hspace{0.45cm}Puma8NH & 0.1340 & 0.1334 & 0.1332 & 0.1334 & 0.1296 & 0.1269 &	0.1295 & 0.1284 \\ \hline
		\hspace{0.7cm}Kinematics & 0.0835 & 0.0831 & 0.0830 & 0.0831 & 0.0804 & 0.0801 &	0.0803 & 0.0801 \\ \hline
		\hspace{0.55cm}Compactiv & 0.0606 & 0.0599 & 0.0608 & 0.0598 & 0.0137 & 0.0086 &	0.0304 & 0.0078 \\ \hline
		\hspace{0.55cm}Protein Tertiary & 0.2554 & 0.2550 & 0.2549 & 0.2550 & 0.2370 & 0.2334 &	0.2385 & 0.2373 \\ \hline
		\hspace{0.55cm}ONP & 0.0015 & 0.0009 & 0.0009 & 0.0009 & 0.0009 & 0.0009 &	0.0009 & 0.0009 \\ \hline
		\hspace{0.55cm}California Housing & 0.0446 & 0.0450 & 0.0452 & 0.0448 & 0.0685 & 0.0671 &	0.0579 & 0.0683 \\ \hline
		\hspace{0.55cm}YPMSD & 0.0237 & 0.0237 & 0.0233 & 0.0237 & 0.0454 & 0.0337 & 0.0302 & 0.0292 \\ \hline
	\end{tabular}
	\caption{\footnotesize{The MSE of the proposed algorithms on real data sets.}} \label{tab:result}
\end{table*} 

\begin{figure*}
	\begin{subfigure}[b]{0.33\textwidth}
		\centering
		\includegraphics[scale=0.45]{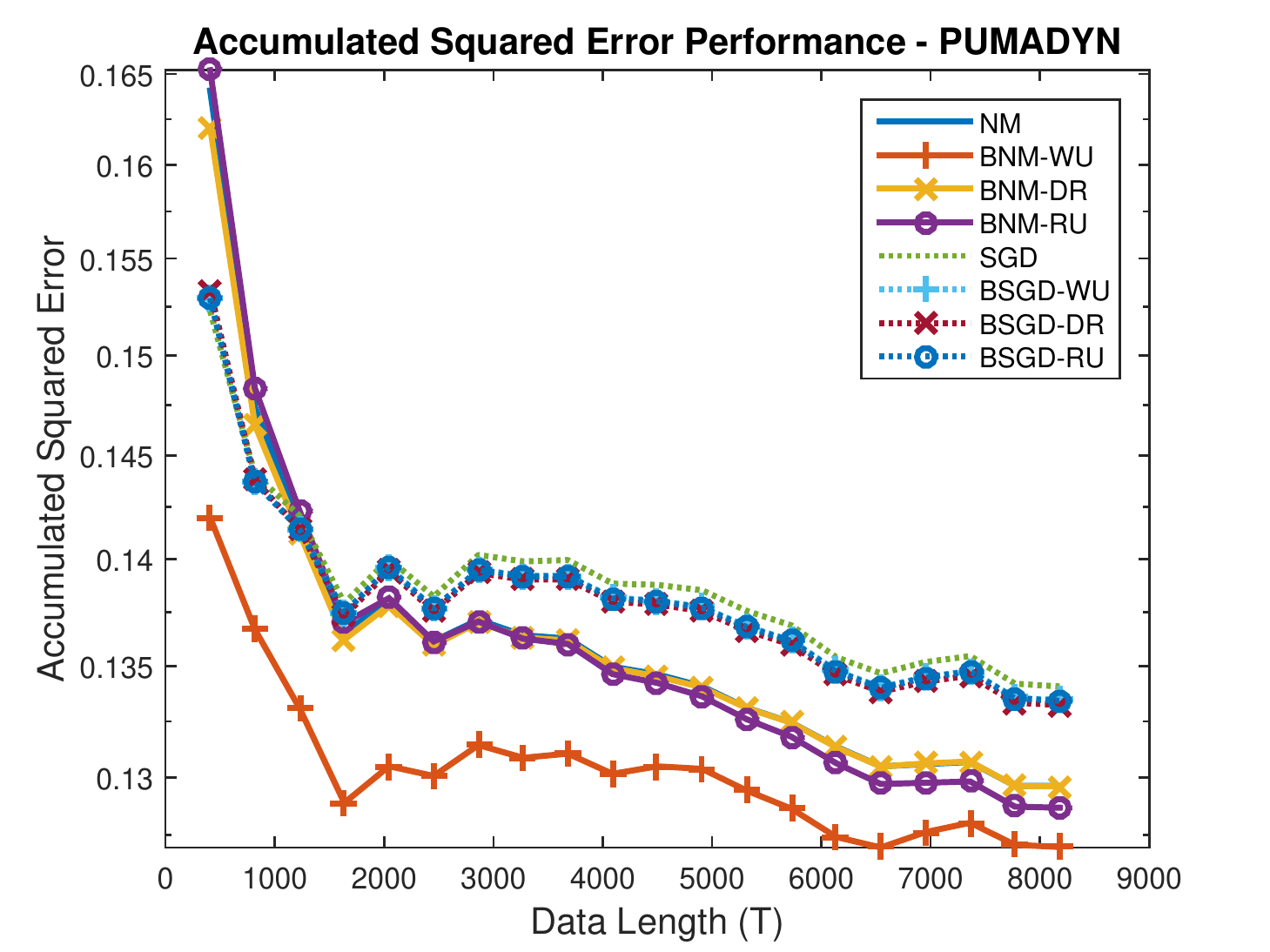}\\
		\caption{Puma8NH}\label{fig:pumadyn}
	\end{subfigure}
	\begin{subfigure}[b]{0.33\textwidth}
		\centering
		\includegraphics[scale=0.45]{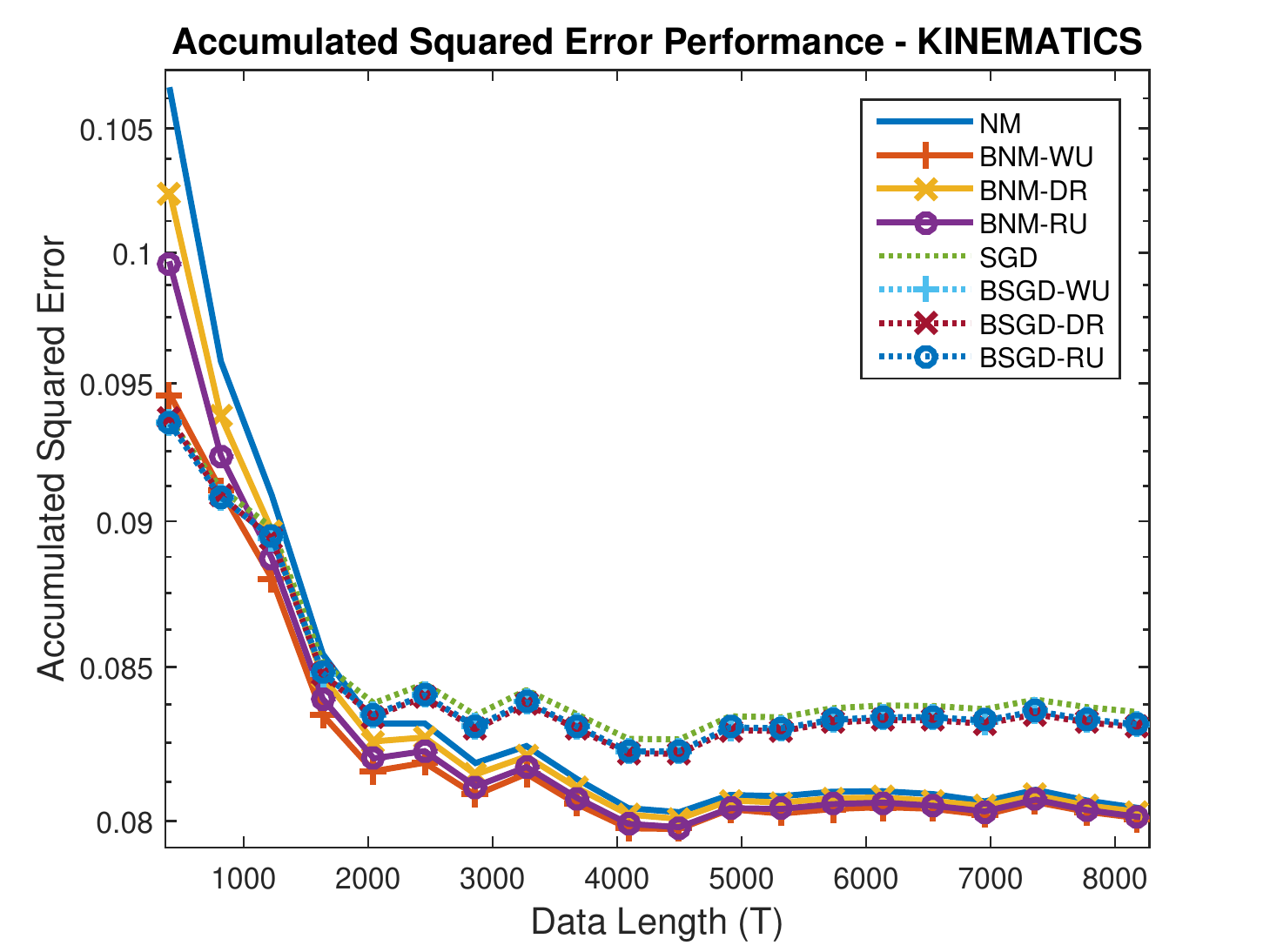}\\
		\caption{Kinematics.}\label{fig:kinematics}
	\end{subfigure}
	\begin{subfigure}[b]{0.33\textwidth}
		\centering
		\includegraphics[scale=0.45]{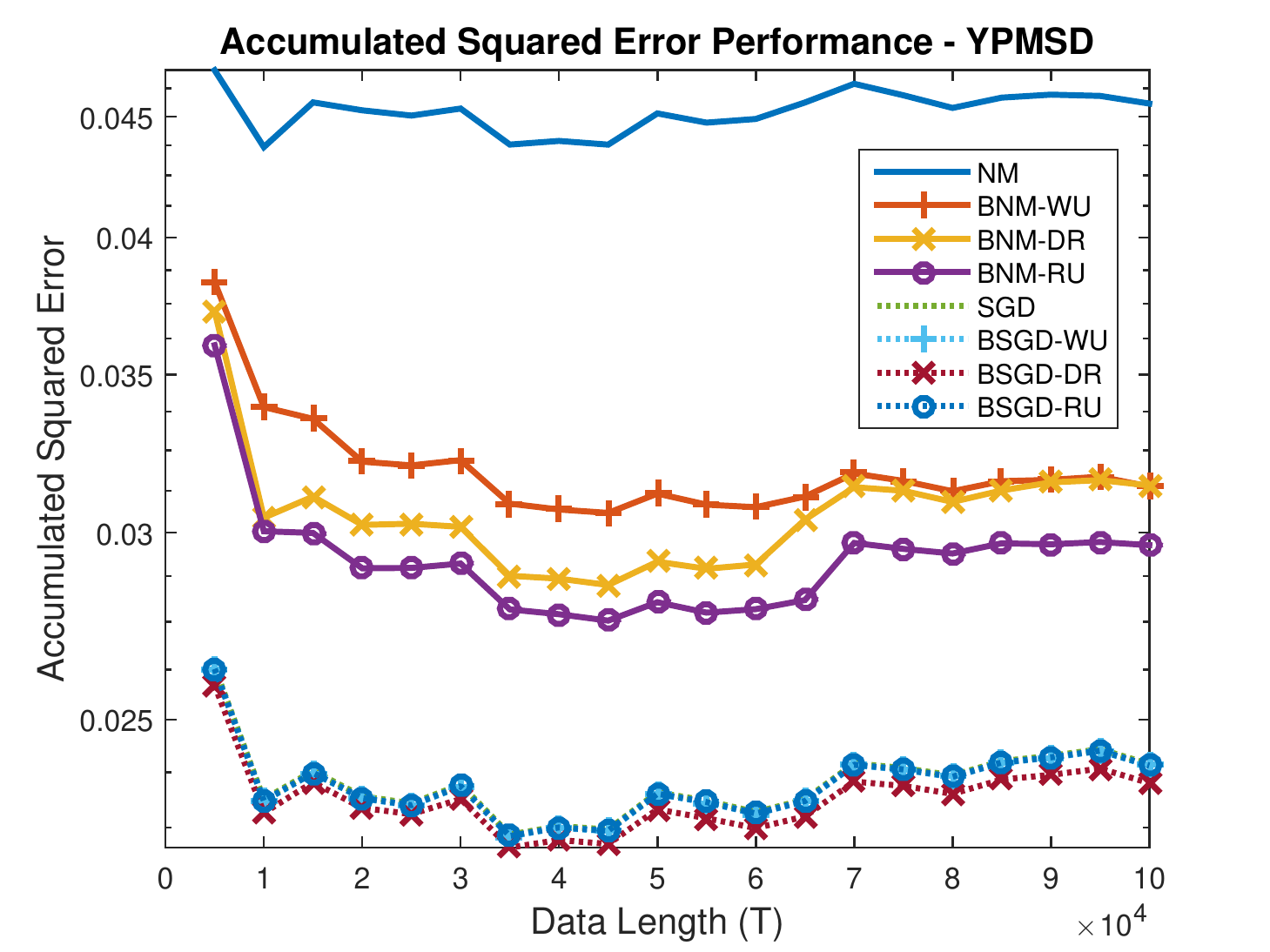}\\
		\caption{YPMSD}\label{fig:ypmsd}
	\end{subfigure}
	\caption{The performance of the proposed boosting methods on three real life data sets.}
\end{figure*}

Here, we briefly explain the details of the data sets:
\begin{enumerate}
	\item MV: This is an artificial dataset with dependencies between the attribute values. One can refer to (\cite{ltorgo}) for further details. There are $10$ attributes and one target value. In this dataset, we can slightly improve the performance of a single linear WL by using any of the proposed methods.
	\item Puma Dynamics (Puma8NH): This dataset is a realistic simulation of the dynamics of a Puma $560$ robot arm (\cite{ltorgo}). The task is to predict the angular acceleration of one of the robot arm's links. The inputs include angular positions, velocities and torques of the robot arm. According to the ASE results in Fig. \ref{fig:pumadyn}, the BNM-WU has the best boosting performance in this experiment. Nonetheless, the SGD-based methods also improve the performance.
	\item Kinematics: This dataset is concerned with the forward kinematics of an 8 link robot arm (\cite{ltorgo}). We use the variant 8nm, which is highly non-linear and noisy. As shown in Fig. \ref{fig:kinematics}, our proposed algorithms slightly improve the performance in this experiment.
	\item Computer Activity (Compactiv): This real dataset is a collection of computer systems activity measures (\cite{ltorgo}). The task is to predict USR, the portion of time that CPUs run in user mode from all attributes (\cite{ltorgo}). The NM-based boosting algorithms deliver a significant performance improvement in this experiment, as shown by the results in Table \ref{tab:result}.
	\item Protein Tertiary (\cite{Lichman2013}): This dataset is collected from Critical Assessment of protein Structure Prediction (CASP) experiments $5-9$. The aim is to predict the size of the residue using $9$ attributes over $45730$ data instances.
	\item Online News Popularity (ONP) (\cite{Lichman2013, pereira2015}): This dataset summarizes a heterogeneous set of features about articles published by Mashable in a period of two years. The goal is to predict the number of shares in social networks (popularity).
	\item California Housing: This dataset has been obtained from StatLib repository. They have collected information on the variables using all the block groups in California from the 1990 Census. Here, we seek to find the house median values, based on the given attributes. For further description one can refer to (\cite{ltorgo}).
	\item Year Prediction Million Song Dataset (YPMSD) (\cite{Bertin-Mahieux2011}): The aim is predicting the release year of a song from its audio features. Songs are mostly western, commercial tracks ranging from 1922 to 2011, with a peak in the year 2000s. We use a subset of the Million Song Dataset (\cite{Bertin-Mahieux2011}). As shown in Table \ref{tab:result} and Fig. \ref{fig:ypmsd}, our algorithms can significantly improve the performance of the linear WL in this experiment.
\end{enumerate}

%\begin{figure*}
%	\centering
%	\begin{subfigure}[b]{0.4\textwidth}
%		\centering
%		\includegraphics[scale=0.5]{figures/new_results/cal_housing_rel.eps}
%		\caption{California Housing}\label{fig:california}
%	\end{subfigure}
%	\begin{subfigure}[b]{0.4\textwidth}
%		\centering
%		\includegraphics[scale=0.5]{figures/new_results/compactiv_lms_rel.eps}
%		\caption{Computer Activity (CompAct)}\label{fig:compact_lms}
%	\end{subfigure}
%	\linebreak
%	\begin{subfigure}[b]{0.4\textwidth}
%		\centering
%		\includegraphics[scale=0.5]{figures/new_results/compactiv_rls_rel.eps}
%		\caption{Computer Activity (CompAct)}\label{fig:compact_rls}
%	\end{subfigure}
%	\begin{subfigure}[b]{0.4\textwidth}
%		\centering
%		\includegraphics[scale=0.5]{figures/new_results/bank8fm_rel.eps}
%		\caption{Bank data}\label{fig:bank}
%	\end{subfigure}
%	\caption{The performance results of the boosting approaches in the benchmark real and synthetic data experiment.}\label{fig:real_results}
%\end{figure*}

%\begin{figure}
%	\centering
%	\includegraphics[width=0.4\textwidth]{figures/new_results/compactiv_rls_rel.eps}\\
%	\caption{The performance of RLS-based boosting methods on the CompAct data set.}\label{fig:compact_rls}
%\end{figure}
%
%\begin{figure}
%	\centering
%	\includegraphics[width=0.4\textwidth]{figures/new_results/bank8fm_rel.eps}\\
%	\caption{The performance of LMS-based boosting methods on the Bank data set.}\label{fig:bank}
%\end{figure}

\section{Conclusion}
We introduced a novel family of boosted online regression algorithms and proposed three different boosting approaches, i.e., weighted updates, data reuse, and random updates, which can be applied to different online learning algorithms. We provide theoretical bounds for the MSE performance of our proposed methods in a strong mathematical sense. We emphasize that while using the proposed techniques, we do not assume any prior information about the statistics of the desired data or feature vectors. We show that by the proposed boosting approaches, we can significantly improve the MSE performance of the conventional SGD and NM algorithms. Moreover, we provide an upper bound for the weights generated during the algorithm that leads us to a thorough analysis of the computational complexity of these methods. The computational complexity of the random updates method is remarkably lower than that of the conventional mixture-of-experts and other variants of the proposed boosting approaches, without degrading the performance. Therefore, the boosting using random updates approach is an elegant alternative to the conventional mixture-of-experts method when dealing with real life large scale problems. We provide several results that demonstrate the strength of the proposed algorithms over a wide variety of synthetic as well as real data.
%We introduce a novel family of boosted online regression algorithms and propose three different boosting approaches, ``weighted updates'',``data reuse'', and ``random updates'', which are applicable to different online learning algorithms. We provide strong mathematical justifications for the MSE performance of our proposed methods, in an individual sequence manner. We show that by the proposed boosting approaches, we can significantly improve the MSE performance of the conventional SGD and NM algorithms. Moreover, we provide an upper bound for the weights generated during the algorithm, which leads us to a thorough analysis of the complexity of these methods. The complexity of random updates method is remarkably lower than that of the mixture-of-experts and other two boosting approaches, while the MSE performance does not degrade. Therefore, the boosting using random updates approach is an elegant alternative to the conventional mixture-of-experts method when dealing with real life large scale problems.\par

\section*{Appendices}
\appendix

\section{Proof of Lemma 1.}\label{app:lem1}
We observe that according to Algorithm \ref{alg:BORA},
\begin{align*}
l_t^{(M+1)} & = \sum_{k=1}^{M} [\sigma_m^2 - (e_t^{(k)})^2],\\
e_t & = \frac{1}{M} \sum_{k=1}^{M} e_t^{(k)},
\end{align*}
%\begin{align*}
%e_t & = d_t - \hat{d}_t \\
%& = \frac{1}{M} \sum_{k=1}^{M} [d_t - \hat{d}_t^{(k)}]\\
%& = \frac{1}{M} \sum_{k=1}^{M} e_t^{(k)},
%\end{align*}
In addition, we have
\[
\sum_{k=1}^{M} \left( e_t^{(k)} \right)^2 \geq \frac{1}{M} \left( \sum_{k=1}^{M} e_t^{(k)} \right)^2,
\]
and as a result, if $e_t^2 > \sigma_m^2$, then $l_t^{(M+1)} \leq 0$, i.e., $\lambda_t^{(M+1)}=1$. Hence by defining a \textit{modified} desired MSE as $\sigma_m^2 \triangleq \frac{\sigma_d^2 - \kappa}{1 - \kappa}$, and $\vz_t = [1/M, ..., 1/M]$ for $t = 1, ...,T$, we have
\begin{align*}
\frac{|\{t:e_t^2 > \sigma_m^2\}|}{T} & = \frac{|\{t:\lambda_t^{(M+1)}=1\}|}{T} \\
& \leq \frac{\sum_{t=1}^{T} \lambda_t^{(M+1)}}{T}\\
& \leq \kappa.
\end{align*}
Finally we have
\begin{align*}
\frac{\sum_{t=1}^{T}e_t^2}{T} & = \frac{\sum_{t; e_t^2 \leq \sigma_m^2}e_t^2}{T} + \frac{\sum_{t; e_t^2 > \sigma_m^2}e_t^2}{T} \\
& \leq \frac{\sum_{t; e_t^2 \leq \sigma_m^2}\sigma_m^2}{T} + \frac{\sum_{t; e_t^2 > \sigma_m^2} 1}{T}\\
& \leq (1-\kappa)\sigma_m^2 + \kappa\\ 
& = \sigma_d^2.
\end{align*}
This completes the proof of Lemma 1. $\Box$

\section{Proof of Lemma 2.}\label{app:lem2}
We have
\[
\sum_{t=1}^{T} \sum_{k=1}^{M} \lambda_t^{(k)} \left[ 1-\left(e_t^{(k)}\right)^2 \right] \geq (1-4\sigma^2) \sum_{t=1}^{T} \sum_{k=1}^{M} \lambda_t^{(k)}.
\]
Moreover, since $0 \geq -\left(e_t^{(k)}\right)^2 = l_t^{(k+1)}-l_t^{(k)}-\sigma_m^2 \geq -4$, following the similar lines as the proof of Lemma 5 in (\cite{smooth_boost}), we find that
\[
\sum_{t=1}^{T} \sum_{k=1}^{M} \lambda_t^{(k)} \left[ 1-\left(e_t^{(k)}\right)^2 \right] \leq -\sigma^4 \sigma_m^2 \sum_{t=1}^{T} \sum_{k=1}^{M} \lambda_t^{(k)} + \frac{1}{\sigma \ln(1/\sigma)}.
\]
Since $\sum_{t=1}^{T} \sum_{k=1}^{M} \lambda_t^{(k)} \geq \kappa T M$, we conclude that
\[
M \leq \frac{1}{(\kappa \sigma \ln(1/\sigma))(1-4\sigma^2+\sigma^4 \sigma_m^2)}.
\]
This concludes the proof of Lemma 2. $\Box$
\section*{Acknowledgments}
This work is supported in part by Turkish Academy of Sciences Outstanding Researcher Programme, TUBITAK Contract No. 113E517, and Turk Telekom Communications Services Incorporated.

%\bibliographystyle{MYplainnat}
%\bibliography{myboosting_references}

%\bibliographystyle{plain}
%\input{manuscript.bbl}

\end{document}